\input amstex
\documentstyle{amsppt}
\document

\magnification 1100

\def\gen{{\frak{g}}}

\def\a{{\alpha}}
\def\l{{\lambda}}
\def\b{{\beta}}
\def\eps{{\varepsilon}}

\def\1b{{\bold 1}}

\def\eb{{\bold e}}
\def\fb{{\bold f}}
\def\hb{{\bold h}}

\def\kb{{\bold k}}

\def\rb{{\bold r}}

\def\xb{{\bold x}}

\def\Kb{{\bold K}}

\def\Rb{{\bold R}}

\def\Ub{{\bold U}}

\def\Wb{{\bold W}}

\def\Prod{{\ts\prod}}
\def\Wedge{{\ts\bigwedge}}

\def\s{{\roman s}}

\def\top{{\text{top}}}

\def\spec{\text{spec}\,}

\def\Hom{\text{Hom}\,}

\def\Det{\text{Det}\,}

\def\simto{\,{\buildrel\sim\over\to}\,}

\def\End{\text{End}\,}

\def\Ker{\text{Ker}\,}

\def\rk{\text{rk\,}}

\def\GL{{\text{GL}}}

\def\AA{{\Bbb A}}

\def\CC{{\Bbb C}}
\def\DD{{\Bbb D}}

\def\NN{{\Bbb N}}

\def\ZZ{{\Bbb Z}}

\def\Bc{{\Cal B}}

\def\Ec{{\Cal E}}
\def\Fc{{\Cal F}}

\def\Lc{{\Cal L}}

\def\Nc{{\Cal N}}
\def\Oc{{\Cal O}}

\def\Tc{{\Cal T}}

\def\Vc{{\Cal V}}
\def\Wc{{\Cal W}}

\def\and{{\quad\text{and}\quad}}

\def\ts{\textstyle}

\def\qed{\hfill $\sqcap \hskip-6.5pt \sqcup$}        
\overfullrule=0pt                                    

\def\u1{{\underline 1}}

\def\oh{{\overline h}}

\def\la{{\langle}}
\def\ra{{\rangle}}

\newdimen\Squaresize\Squaresize=14pt
\newdimen\Thickness\Thickness=0.5pt
\def\Square#1{\hbox{\vrule width\Thickness
	      \alphaox to \Squaresize{\hrule height \Thickness\vss
	      \hbox to \Squaresize{\hss#1\hss}
	      \vss\hrule height\Thickness}
	      \unskip\vrule width \Thickness}
	      \kern-\Thickness}
\def\Vsquare#1{\alphaox{\Square{$#1$}}\kern-\Thickness}

\title Standard modules of quantum affine algebras\endtitle
\rightheadtext{standard modules}
\author M. Varagnolo and E. Vasserot\endauthor
\thanks
Both authors are partially supported by EEC grant
no. ERB FMRX-CT97-0100.\endthanks
\abstract
We give a proof of the cyclicity conjecture of Akasaka-Kashiwara,
for simply laced types, via quiver varieties. We get also 
an algebraic characterization of the standard modules.
\endabstract
\endtopmatter
\document

\head 1. Introduction\endhead
Let $\gen$ be a simple, simply laced, complex Lie algebra.
If $\gen$ is of type $A$, a geometric realization
of the quantized enveloping algebra $\Ub$ of $\gen[t,t^{-1}]$ and of its 
simple modules was given a few years ago in \cite{GV}, \cite{V}. 
This construction involved perverse sheaves and the convolution algebra in 
equivariant $K$-theory of partial flags varieties of type $A$.
It was then observed in \cite{N1}, \cite{N2}, that these varieties should
be viewed as a particular case of the quiver varieties
associated to any symmetric Kac-Moody Lie algebra. This leads to a geometric
realization of $\Ub$ via a convolution algebra in equivariant
$K$-theory of the quiver varieties for $\gen$ of (affine) type $A^{(1)}$ 
in \cite{VV} and for a general symmetric Kac-Moody algebra $\gen$ in \cite{N3}.
For any symmetric Lie algebra $\gen$, one gets a formula for the 
dimension of the finite dimensional 
simple modules of $\Ub$ in terms of intersection cohomology
(see \cite{N3}). 
A basic tool in this geometric approach are the standard modules.
They are the geometric counterpart of the Weyl modules of $\Ub$
(see Remark 7.19).
In this paper we give an algebraic construction of the standard modules. 
It answers to a question in \cite{N3} (Corollary 7.16). An immediate corollary
is a proof of the cyclicity conjecture in \cite{AK} (Corollary 7.17) for simply
laced Lie algebras.

The plan of the paper is the following : Sections 1 to 6 contain recollections
on quantum affine algebras and quiver varieties. 
The main results are given in Section 7. 
The proof of Theorem 7.4 uses Lemma 8.1.

While we were preparing this paper Kashiwara mentioned to us that he
has proved the conjecture in \cite{AK} by a different approach
(via canonical bases). We would like to thank the referee for numerous 
remarks on the first version of the paper.

\vskip3mm

\head 2. The algebra $\Ub$\endhead
Let $\gen$ be a simple, simply laced, complex Lie algebra.
The quantum loop algebra associated to ${\frak g}$ is the
$\CC(q)$-algebra $\Ub'$ generated by 
$\xb^\pm_{ir},\,\kb^\pm_{is},\,\kb_i^{\pm 1}=\kb^\pm_{i0}$ 
$(i\in I,\,r\in\ZZ,\,s\in\pm\NN^\times)$
modulo the following defining relations 
$$\kb_i\kb_i^{-1}=1=\kb^{-1}_i\kb_i,
\quad [\kb^\pm_{i,\pm r},\kb^\eps_{j,\eps s}]=0,$$
$$\kb_i\xb^\pm_{jr}\kb_i^{-1}=q^{\pm a_{ij}}\xb^\pm_{jr},$$
$$(w-q^{\pm a_{ji}}z)\,\kb^\eps_j(w)\,\xb^\pm_i(z)=
(q^{\pm a_{ji}}w-z)\,\xb^\pm_i(z)\,\kb^\eps_j(w),$$
$$(z-q^{\pm a_{ij}}w)\xb^\pm_i(z)\xb^\pm_j(w)=
(q^{\pm a_{ij}}z-w)\xb^\pm_j(w)\xb^\pm_i(z),$$
$$[\xb^+_{ir},\xb^-_{js}]=\delta_{ij}{{\kb_{i,r+s}^+-\kb_{i,r+s}^-}
\over q-q^{-1}},$$
$$\sum_w\sum_{p=0}^m(-1)^{^p}
\left[\matrix m\cr p\endmatrix\right]
\xb^\pm_{ir_{w(1)}}\xb^\pm_{ir_{w(2)}}\cdots\xb^\pm_
{ir_{w(p)}}\xb^\pm_{js}\xb^\pm_{ir_{w(p+1)}}\cdots\xb^\pm_{ir_{w(m)}}=0,$$
where $i\neq j,$ $m=1-a_{ij},$ $r_1,...,r_m\in\ZZ,$ and $w\in S_m$.
We have set $[n]=q^{1-n}+q^{3-n}+...+q^{n-1}$ if $n\geq 0$,
$[n]!=[n][n-1]...[2]$, and
$$\left[\matrix m\cr p\endmatrix\right]=
{[m]!\over[p]![m-p]!}.$$
We have also set $\eps=+$ or $-$,
$$\kb^\pm_i(z)=\sum_{r\geq 0}\kb^\pm_{i,\pm r}z^{\mp r},\quad
\xb_i^\pm(z)=\sum_{r\in\ZZ}\xb_{ir}^\pm\,z^{\mp r}.$$
Put $\AA=\CC[q,q^{-1}]$.
Consider also the $\AA$-subalgebra $\Ub\subset\Ub'$ 
generated by the quantum divided powers
$\xb^{\pm\,(n)}_{ir}=\xb^{\pm\,n}_{ir}/[n]!$, the Cartan elements
$\kb_i^{\pm 1}$, and the elements $\hb_{is}$ such that 
$$\kb^\pm_i(z)=\kb_i^{\pm 1}\exp
\Bigl(\pm(q-q^{-1})\sum_{s\geq 1}\hb_{i,\pm s}z^{\mp s}\Bigr).$$
Let $\Delta^\circ$ be the coproduct defined in terms of the Kac-Moody generators
$\eb_i,\fb_i,\kb_i^{\pm 1}$, $i\in I\cup\{0\}$, of $\Ub'$ as follows
$$\Delta^\circ(\eb_i)=\eb_i\boxtimes 1+\kb_i\boxtimes\eb_i,\quad 
\Delta^\circ(\fb_i)=\fb_i\boxtimes\kb_i^{-1}+1\boxtimes\fb_i,\quad
\Delta^\circ(\kb_i)=\kb_i\boxtimes\kb_i,$$
where $\boxtimes$ is the tensor product over the field $\CC$, or the ring $\AA$.
Let $\tau$ be the anti-automorphism of $\Ub$ such that
$\tau(\eb_i)=\fb_i$, $\tau(\fb_i)=\eb_i$, $\tau(\kb_i)=\kb^{-1}_i$,
and $\tau(q)=q^{-1}$. It is known (see \cite{B}) that
$\tau(\xb^-_{i,-k})=\xb^+_{ik}$ and $\tau(\kb^\pm_{i,\pm r})=\kb^\mp_{i,\mp r}$.
Let $\Delta^\bullet$ be the coproduct opposit to $\Delta^\circ$. 
We have $(\tau\boxtimes\tau)\Delta^\circ\tau=\Delta^\bullet$.
Hereafter $\zeta$ is an element of $\CC^\times$ which is not a root of unity. 

\vskip3mm

\head 3. The braid group\endhead
Let $W, P,Q,$ be the Weyl group, the weight lattice, 
and the root lattice of $\gen$.
The extended affine Weyl group $\tilde W=W\ltimes P$
is generated by the simple reflexions
$s_i$ and the fundamental weights $\omega_i$, $i\in I$.
Let $\Gamma$ be the quotient of $\tilde W$ by the normal Coxeter
subgroup generated by the simple affine reflexions.
The group $\Gamma$ is a group of diagram automorphisms of the 
extended Dynkin diagram of $\gen$. In particular $\Gamma$ acts 
on $\Ub$ and on $\tilde W$.
For any $w\in\tilde W$ let $l(w)$ denote its length.
The braid group $B_{\tilde W}$ associated to $\tilde W$ is the 
group on generators $T_w,$ $w\in\tilde W$, 
with the relation $T_wT_{w'}=T_{ww'}$ whenever $l(ww')=l(w)+l(w')$.
The group $B_{\tilde W}$ acts on $\Ub$ by algebra 
automorphisms (see \cite{L1}, \cite{B}). 
Recall that $\xb^-_{ir}=\nu(i)^rT_{\omega_i}^r(\fb_i)$, 
$\xb^+_{ir}=\nu(i)^rT_{\omega_i}^{-r}(\eb_i)$,
for a fixed function $\nu\,:\,I\to\{\pm 1\}$ such that $\nu(i)+\nu(j)=0$
if $a_{ij}<0$ (see \cite{B, Definition 4.6}).
We have $T_{\omega_j}(\kb^\pm_{ir})=\kb^\pm_{ir}$ for any $i,j$,
and $T_{\omega_j}(\xb_{ir}^\pm)=\xb_{ir}^\pm$ for any $i\neq j$.
For any $i\in I\cup\{0\}$ put
$$R_i=\sum_{l\geq 0}c_l\,T_{s_i}(\fb_i)^{(l)}\boxtimes T_{s_i}(\eb_i)^{(l)}\quad\text{where}\quad
c_l=(-1)^lq^{-l(l-1)/2}(q-q^{-1})^l[l]!.$$
The element $R_i$ belongs to a completed tensor product
$(\Ub\boxtimes\Ub)\hat{}$ (see \cite{L1, \S 4.1.1} for instance).
It is known that
$$R_i^{-1}=\sum_{l\geq 0}\bar c_l\,T_{s_i}(\fb_i)^{(l)}\boxtimes T_{s_i}(\eb_i)^{(l)}\quad\text{where}\quad
\bar c_l=q^{l(l-1)/2}(q-q^{-1})^l[l]!.$$
If $\tau s_{i_1}...s_{i_r}$ is a reduced expression of $w\in\tilde W$
with $\tau\in\Gamma$, set
$$R_w=\tau\bigl(T^{[2]}_{i_1}...T^{[2]}_{i_{r-1}}(R_{i_r})...
T^{[2]}_{i_1}(R_{i_2})R_{i_1}\bigr),\leqno(3.1)$$
where $a^{[2]}=a\boxtimes a$ for any $a$.
In particular if $w,w'\in\tilde W$ are such that $l(ww')=l(w)+l(w')$, then
$R_{ww'}=T^{[2]}_w(R_{w'})R_w.$
For any $w\in\tilde W$ set $\Delta^\circ_w=T^{[2]}_w\Delta^\circ T^{-1}_w.$
Then
$R_w\cdot\Delta^\circ x\cdot R_w^{-1}=\Delta^\circ_wx$, for all $x\in\Ub$
(see \cite{L1, 37.3.2} and \cite{B, Section 5}).

\vskip3mm

\head 4. The quiver varieties\endhead
\subhead 4.1\endsubhead
Let $I$ (resp. $E$) be the set of vertices (resp. edges) of a finite
graph $(I,E)$ with no edge loops. For $i,j\in I$ let $n_{ij}$ be
the number of edges joining $i$ and $j$. Put $a_{ij}=2\delta_{ij}-n_{ij}$.
The map $(I,E)\mapsto A=(a_{ij})_{i,j\in I}$ is a bijection from the set of
finite graphs with no edge loops onto the set of symmetric generalized 
Cartan matrices. Let $\a_i$, $i\in I$,
be the simple roots of the symmetric Kac-Moody 
algebra $\gen$ corresponding to $A.$
Hereafter we assume that $\gen$ is finite dimensional,
i.e. the matrix $A$ is positive definite. 
Let $H$ be the set of edges of $(I,E)$ together 
with an orientation. For $h\in H$ let $h'$ and $h''$
the incoming and the outcoming vertex of $h$. If $h\in H$ we denote
by $\oh\in H$ the same edge with opposite orientation.
Given two $I$-graded finite dimensional complex vector spaces 
$V=\bigoplus_{i\in I}V_i,\, W=\bigoplus_{i\in I}W_i$, set
$$E(V,W)=\bigoplus_{h\in H}\Hom(V_{h''},W_{h'}),\quad
L(V,W)=\bigoplus_{i\in I}\Hom(V_i,W_i).$$
Let $P^+, Q^+$ be the semi-groups
$Q^+=\bigoplus_{i\in I}\NN\a_i$ and $P^+=\bigoplus_{i\in I}\NN\omega_i$.
Let us fix once for all the following convention :
the dimension of the graded vector space $V$
is identified with the element $\alpha=\sum_{i\in I}v_i\a_i\in Q^+$ 
(where $v_i$ is the dimension of $V_i$), 
while the dimension of $W$ is 
identified with the weight $\lambda=\sum_iw_i\omega_i\in P^+$ 
(where $w_i$ is the dimension of $W_i$). 
We also put $|\lambda|=\sum_iw_i$ and $|\alpha|=\sum_iv_i$. 
We write $\alpha\geq\alpha'$ if and only if $\alpha-\alpha'\in Q^+$.
Set
$$M_{\alpha\lambda}=E(V,V)\oplus L(W,V)\oplus L(V,W).$$
For any $(B,p,q)\in M_{\alpha\lambda}$ let $B_h$ be the component of $B$
in $\Hom(V_{h''},V_{h'})$ and set
$$m_{\alpha\lambda}(B,p,q)=\sum_h\eps(h)B_hB_\oh +pq\in L(V,V),$$
where $\eps$ is a function $\eps\,:\, H\to\CC^\times$ such
that $\eps(h)+ \eps(\oh)=0.$ 
A triple $(B,p,q)\in m_{\alpha\lambda}^{-1}(0)$ is stable
if there is no nontrivial $B$-invariant subspace of $\Ker q$.
Let $m_{\alpha\lambda}^{-1}(0)^\circ$ be the set of stable triples. 
The group $G_\alpha=\prod_i\GL(V_i)$ acts on $M_{\alpha\lambda}$ by 
$$g\cdot(B,p,q)=(gBg^{-1},gp,qg^{-1}).$$
The action of $G_\alpha$ on $m_{\alpha\lambda}^{-1}(0)^\circ$ is free.
Put 
$$Q_{\alpha\lambda}=m_{\alpha\lambda}^{-1}(0)^\circ/G_\alpha\and
N_{\alpha\lambda}=m^{-1}_{\alpha\lambda}(0)/\!\!/G_\alpha,$$
where $/\!\!/$ denotes the categorical quotient. 
The variety $Q_{\alpha\lambda}$ is smooth and quasi-projective.

\subhead 4.2\endsubhead
Let $\pi\,:\, Q_{\a\l}\to N_{\a\l}$ be the affinization map. 
Let $L_{\a\l}=\pi^{-1}(0)\subset Q_{\a\l}$ be the zero fiber.
It is known that $\dim Q_{\a\l}=2\dim L_{\a\l}$.
If $\a,\a'\in Q^+$ are such that $\a\geq\a'$, 
then the extension by zero of representations of the quiver
gives an injection $N_{\a'\l}\hookrightarrow N_{\a\l}$
(see \cite{N3, Lemma 2.5.3}). For any $\a,\a',$ we consider the fiber product
$$Z_{\a\a'\l}=Q_{\a\l}\times_\pi Q_{\a'\l}.$$
It is known that $\dim Z_{\a\a'\l}=(\dim Q_{\a\l}+\dim Q_{\a'\l})/2$.
If $\a'=\a+\a_i$ and $V\subset V'$ have dimension $\a$, $\a'$, respectively,
let $C^{i+}_{\a'\l}\subset Z_{\a\a'\l}$
be the set of pairs of triples $(B,p,q),$ $(B',p',q'),$ such that
$B'_{|V}=B,$ $p'=p,$ $q'_{|V}=q.$ 
If $\a'=\a-\a_i$, put
$C^{i-}_{\a'\l}=\phi\bigl(C^{i+}_{\a\l}\bigr)\subset Z_{\a\a'\l}$,
where $\phi$ flips the components.
The variety $C^{i\pm}_{\a'\l}$ is an
irreducible component of $Z_{\a\a'\l}$. 
Consider the following varieties 
$$N_\l=\bigcup_\a N_{\a\l},\quad Q_\l=\bigsqcup_\a Q_{\a\l},\quad
Z_\l=\bigsqcup_{\a,\a'}Z_{\a\a'\l},\quad
C^{i\pm}_\l=\bigsqcup_\a C^{i\pm}_{\a\l},\quad L_\lambda=\bigsqcup_\a L_{\a\l},$$
where $\a,\a',$ take all the possible values in $Q^+$.
Observe that for a fixed $\l$, the set $Q_{\a\l}$ is empty except
for a finite number of $\a$'s.

\subhead 4.3\endsubhead
Put $\tilde{G}_\l=G_\l\times\CC^\times.$
The group $\tilde G_\l$ acts on $M_{\a\l}$ by 
$$(g,z)\cdot (B,p,q)=(z B,z pg^{-1},z gq).$$
This action descends to $Q_{\a\l}$ and $N_{\a\l}$.
For any element $s=(t,\zeta)\in\tilde G_\l$ with $t$ semi-simple, 
let $\la s\ra\subset\tilde G_\l$
be the Zariski closure of $s^\ZZ$.
For any group homomorphism $\rho\in\Hom(\la s\ra,G_\a)$, let 
$Q(\rho)\subset Q_{\a\l}$ be the subset of the classes of
the triples $(B,p,q)$ such that 
$s\cdot (B,p,q)=\rho(s)\cdot (B,p,q).$
The fixpoint set $Q_{\a\l}^s$ is the disjoint union of the subvarieties
$Q(\rho)$. It is proved in \cite{N3, Theorem 5.5.6} that 
$Q(\rho)$ is either empty or a connected component of $Q_{\a\l}^s$.

\proclaim{Lemma 4.4}
(i) Fix $\lambda^1,\lambda^2\in P^+,$ such that $\lambda=\lambda^1+\lambda^2$. 
The direct sum $M_{\l^1}\times M_{\l^2}\to M_\l$ gives a closed embedding
$\kappa\,:\,Q_{\l^1}\times Q_{\l^2}\hookrightarrow Q_\l.$\hfill\break
(ii) Fix a semi-simple element
$t=t^1\oplus t^2\in G_{\l^1}\times G_{\l^2}$. 
Set $s=(t,\zeta)\in\tilde G_\l$, 
$s^1=(t,\zeta)\in\tilde G_{\l^1}$, 
$s^2=(t,\zeta)\in\tilde G_{\l^2}$. 
If $(\zeta^\ZZ\spec t^1)\cap\spec t^2=\emptyset$, then 
$Q_\l^s=\kappa(Q^{s^1}_{\l^1}\times Q^{s^2}_{\l^2})$. 
\endproclaim

\noindent{\sl Proof of 4.4.}
Fix $I$-graded vector spaces $W,W^1,W^2,V,V^1,V^2,$ such that
$W=W^1\oplus W^2$, $V=V^1\oplus V^2$, $\dim W^1=\l^1$, 
$\dim W^2=\l^2$, $\dim V^1=\a^1$, and $\dim V^2=\a^2$.
Fix triples $x^1, {x'}^1\in m_{\a^1\l^1}^{-1}(0)^\circ,$
$x^2, {x'}^2\in m_{\a^2\l^2}^{-1}(0)^\circ$. 
The triple $x=x^1\oplus x^2$ is stable since if $V'\subseteq\Ker q$ is
a $B$-stable subspace, then $V^1\cap V'=\{0\}$ by the stability of $x^1$,
and then $V'$ embeds in $V/V^1$. Thus it is zero by the stability of $x^2$. 
Assume that $g\in G_{\a^1+\a^2}$ maps $x$ to ${x'}^1\oplus{x'}^2.$ 
Then,
$$q(g^{-1}(V^2)\cap V^1)\subseteq W^1\cap W^2=\{0\}\and
B(g^{-1}(V^2)\cap V^1)\subseteq g^{-1}(V^2)\cap V^1.$$
Thus the stability of $x$ gives $g^{-1}(V^2)\cap V^1=\{0\}.$
In the same way we get $g^{-1}(V^1)\cap V^2=\{0\}.$
Thus $g\in G_{\a^1}\times G_{\a^2}$. 
Claim (i) is proved. 

Fix $\rho$ such that $Q(\rho)$ is non empty.
For any $z\in\CC^\times$ put 
$$V(z)=\Ker(\rho(s)-z^{-1} id_V),\qquad W(z)=\Ker(t-z\, id_W).$$
If  $x\in Q(\rho)$ and $(B,p,q)$ is a representative of $x$, then
$$B(V(z))\subset V(z/\zeta),\quad q(V(z))\subset W(z/\zeta),\quad p(W(z))
\subset V(z/\zeta).$$
Moreover the stability of $(B,p,q)$ implies that 
$V=\bigoplus_{z\in\zeta^\ZZ \spec t}V(z).$
Claim (ii) follows.\qed

\vskip3mm

\subhead 4.5\endsubhead
Let $\Vc=m_{\a\l}^{-1}(0)^\circ\times_{G_\a}V$ and 
$\Wc$ be respectively the tautological bundle 
and the trivial $W$-bundle on $Q_{\a\l}.$ 
The $i$-th component of $\Vc,\Wc,$ is denoted by $\Vc_i,\Wc_i$.
The bundles $\Vc,\Wc,$ are $\tilde G_\l$-equivariant. 
Let $q$ be the trivial line bundle on $Q_{\a\l}$
with the degree one action of $\CC^\times$. We consider the classes 
$$\matrix
\Fc^i_{\a\l}&=q^{-1}\Wc_i-(1+q^{-2})\Vc_i+q^{-1}\sum_{h'=i}\Vc_{h''}
\hfill\cr
\Tc_{\a\l}&=qE(\Vc,\Vc)+q^2L(q\Wc-\Vc,\Vc)+L(\Vc,q\Wc-\Vc)\hfill
\endmatrix$$
in $\Kb^{\tilde G_\l}(Q_{\a\l})$, and the classes 
$$\matrix
\Tc^{i+}_{\a'\l}=qE(\Vc,\Vc')+q^2L(q\Wc-\Vc,\Vc')+L(\Vc,q\Wc-\Vc')-q^2
\hfill\cr
\Tc^{i-}_{\a'\l}=qE(\Vc',\Vc)+q^2L(q\Wc-\Vc',\Vc)+L(\Vc',q\Wc-\Vc)-q^2
\hfill
\endmatrix$$
in $\Kb^{\tilde G_\l}(Q_{\a\l}\times Q_{\a'\l})$ (where $\a'=\a\pm\a_i$).
It is known that $\Tc_{\a\l}$ is the class of the tangent sheaf to 
$Q_{\a\l}$, and that $\Tc^{i\pm}_{\a'\l}|_{C^{i\pm}_{\a'\l}}$ 
is the class of the normal sheaf of $C^{i\pm}_{\a'\l}$ 
in $Q_\l\times Q_\l$ (see \cite{N2}).

\vskip3mm

\head 5. The convolution product\endhead
\subhead 5.1\endsubhead
For any complex algebraic linear group $G$, and
any quasi-projective $G$-variety $X$ let $\Kb^G(X)$, $\Kb_G(X)$
be the complexified Grothendieck groups of $G$-equivariant 
coherent sheaves, and locally free sheaves respectively, on $X$. 
We put $\Rb(G)=\Kb^G(point)$.
For simplicity let $f_*, f^*,\otimes,$ denote the derived functors 
$Rf_*, Lf^*,\otimes^L$ (where $\otimes$ is the tensor product of sheaves
of $\Oc_X$-modules). 
We use the same notation for a sheaf 
and its class in the Grothendieck group. 
Hereafter, elements of $\Kb_G(X)$ may be 
identified with their image in $\Kb^G(X).$
The class of the structural sheaf is simply denoted by 1.
Given smooth quasi-projective $G$-varieties $X_1,X_2,X_3,$
consider the projection $p_{ab}\,:\,X_1\times X_2\times X_3\to X_a\times X_b$
for all $1\leq a<b\leq 3$. Consider closed subvarieties 
$Z_{ab}\subset X_a\times X_b$
such that the restriction of $p_{13}$ to 
$p_{12}^{-1}Z_{12}\cap p_{23}^{-1}Z_{23}$ is proper and maps to $Z_{13}$. 
The convolution product is the map
$$\star\,:\quad\Kb^G(Z_{12})\boxtimes\Kb^G(Z_{23})\to\Kb^G(Z_{13}),\quad
\Ec\boxtimes\Fc\mapsto p_{13\,*}\bigl((p_{12}^*\Ec)\otimes(p_{23}^*\Fc)\bigr).$$
See \cite{CG} for more details.
The flip $\phi\,:\,X_a\times X_b\to X_b\times X_a$ gives a map
$\phi_*\,:\,\Kb^G(Z_{ab})\to\Kb^G(Z_{ba})$.
This maps anti-commutes with $\star$, i.e.
$$\phi_*(x_{12}\star x_{23})=
\phi_*(x_{23})\star\phi_*(x_{12}),\qquad\forall x_{12}, x_{23}.$$

\subhead 5.2\endsubhead
Let $\DD_X$ be the Serre-Grothendieck duality operator on $\Kb^G(X)$
(see \cite{L2, Section 6.10} for instance).
Assume that $X$ is the disjoint union of smooth connected
subvarieties $X^{(i)}$.
Assume also that we have fixed a particular invertible element $q\in\Kb^G(X)$.
Put $d_{X^{(i)}}=\dim X^{(i)}$, $D_{X^{(i)}}=q^{d_{X^{(i)}}}\DD_{X^{(i)}}$,
$D_X=\sum_iD_{X^{(i)}}$.
Let $\Omega_X$ be the determinant of the cotangent bundle to $X$.
Then $\DD_X(\Ec)=(-1)^{\dim X}\Ec^*\otimes\Omega_X$
for any $G$-equivariant locally free sheaf $\Ec$.
Put $Z_{ab}^{(ij)}=Z_{ab}\cap(X_a^{(i)}\times X_b^{(j)})$.
Assume that $\Omega_{X^{(i)}_a}=q^{-d_{X^{(i)}_a}}$ for all $a,i$.  
The operators 
$$D_{Z_{ab}}=\sum_{i,j}q^{d_{ab}^{(ij)}}\DD_{Z^{ij}_{ab}},$$
where $d_{ab}^{(ij)}=\bigl(d_{X^{(i)}_a}+d_{X^{(j)}_b}\bigr)/2,$
are compatible with the convolution product $\star$, i.e.
$$D_{Z_{12}}(x_{12})\star D_{Z_{23}}(x_{23})=
D_{Z_{13}}(x_{12}\star x_{23}),\qquad\forall x_{12}, x_{23}$$
(see \cite{L2, Lemma 9.5} for more details).

\subhead 5.3\endsubhead
If $\Ec$ is a $G$-bundle on $X$, we have the element 
$\Wedge_z(\Ec)=\sum_{i=0}^{\rk\Ec}(\Wedge^i\Ec)\cdot z^i\in\Kb_G(X)[z]$,
where $\Wedge^i\Ec$ is the $i$-th wedge product. Clearly 
$\Wedge_z(\Ec+\Fc)=\Wedge_z(\Ec)\otimes\Wedge_z(\Fc)$ for any $\Ec,\Fc$, and
$$\Wedge_{-1}(\Ec)=(-1)^{\rk\Ec}\Det(\Ec)\otimes\Wedge_{-1}(\Ec^*).\leqno(5.4)$$
Observe also that $\Wedge_z(\Ec)$ admits an inverse in $\Kb_G(X)[[z]]$.
Let $\bar\Rb(G)$ be the fraction field of $\Rb(G)$ and set
$\bar\Kb_G(X)=\Kb_G(X)\boxtimes_{\Rb(G)}\bar\Rb(G)$,
$\bar\Kb^G(X)=\Kb^G(X)\boxtimes_{\Rb(G)}\bar\Rb(G)$.
Assume that $G$ is a diagonalizable group 
and that the set of $G$-fixed points of the restriction $\Ec|_{X^G}$ is $X^G$.
Then $\Wedge_{-1}(\Ec)$ is invertible in $\bar\Kb_G(X)$
by the localization theorem and \cite{CG, Proposition 5.10.3}.
In particular the element
$\Wedge_{-1}(\Fc-\Ec)=\Wedge_{-1}(\Fc)\otimes\Wedge_{-1}(\Ec)^{-1}$ 
is well-defined in $\bar\Kb_G(X)$ for any $\Fc$. 
If $G$ is a product of $\GL_n$'s and $H\subset G$ is a torus,
then $\Wedge_{-1}(\Ec)$ is still invertible in $\bar\Kb_G(X)$
if the set of $H$-fixed points of $\Ec|_{X^H}$ is $X^H$
(use \cite{CG, Theorem 6.1.22} which holds, 
although the group $G$ is not simply-connected).
In the sequel we may identify a $G$-bundle and 
its class in the Grothendieck group.

\subhead 5.5\endsubhead
Assume that the group $G$ is Abelian. 
For any $\Rb(G)$-module $M$ and any $s\in G$,
let $M_s$ be the specialisation of $M$ at the maximal ideal
in $\Rb(G)$ associated to $s$. The localization theorem
gives isomorphisms of modules
$$\iota_*\,:\,\Kb^G(X^s)_s\to\Kb^G(X)_s,\quad
\iota_*\,:\,\bar\Kb^G(X^s)\to\bar\Kb^G(X),$$
where $\iota\,:\,X^s\to X$ is the closed embedding.

\vskip3mm
 
\head 6. Nakajima's theorem\endhead
\subhead 6.1\endsubhead
We fix a subset $H^+\subset H$ such that
$H^+\cap\bar H^+=\emptyset$ and $H^+\cup\bar H^+=H$.
For any $i,j\in I$ let $n_{ij}^+$ be the number of 
arrows in $H^+$ from $i$ to $j$. Put $n^-_{ij}=n_{ij}-n^+_{ij}$.
Observe that $n_{ij}^+=n_{ji}^-$. Put 
$$\Fc^{i-}_{\alpha\lambda}=-\Vc_i+q^{-1}\sum_jn_{ij}^-\Vc_j,\quad
\Fc^{i+}_{\alpha\lambda}=q^{-1}\Wc_i-q^{-2}\Vc_i+q^{-1}\sum_jn_{ij}^+\Vc_j.$$
Let $(\ |\ )\,:\, Q\times P\to\ZZ$ be the pairing such that
$(\a_i|\omega_j)=\delta_{ij}$ for all $i,j\in I$.
The rank of $\Fc^i_{\a\l}$ is $(\a_i|\l-\a)$.
Put $\Fc^{i\pm}_\l=\bigoplus_\a\Fc^{i\pm}_{\a\l}$ and
$\Fc^i_\l=\bigoplus_\a\Fc^i_{\a\l}.$ 
Let $f^i_\l, f^{i\pm}_\l,$ be the diagonal operators 
acting on $\Kb^{\tilde G_\l}(Q_{\a\l})$ by the scalars
$f^i_{\a\l}=\rk\Fc^i_{\a\l}$ and $f^{i\pm}_{\a\l}=\rk\Fc^{i\pm}_{\a\l}$.
Let $p,p'\,:\,(Q_\l)^2\to Q_\l$ be the first and the second projection.
We denote by $\Vc,\Vc'\in\Kb^{\tilde G_\l}\bigl((Q_\l)^2\bigr)$
the pull-back of the tautological sheaf 
(i.e. $\Vc=p^*\Vc$ and $\Vc'={p'}^*\Vc$). 
Set $\Lc=q^{-1}(\Vc'-\Vc)$. For any $r\in\ZZ$ set 
$$x_{ir}^\pm=\sum_{\a'}(\pm\Lc)_{|C^{i\pm}_{\a'\l}}^{\otimes r+f^{i\pm}_{\a'\l}}
\star\delta_*x^{i\pm}_{\a'\l},\quad
k_i^\pm(z)=\delta_*q^{f^i_\l}\Wedge_{-1/z}\bigl((q^{-1}-q)\Fc^i_{\l}\bigr)^\pm,
\leqno(6.2)$$
where 
$x^{i\pm}_{\a\l}=(-1)^{f^{i\pm}_{\a\l}}
\Det(\Fc^{i\pm *}_{\a\l}),$
the map $\delta$ is the diagonal embedding
$Q_\lambda\hookrightarrow(Q_\l)^2$,
and $\pm$ is the expansion at $z=\infty$ or $0$.
Hereafter we may omit $\delta$, 
hoping that it makes no confusion.
Let $\Ub_\lambda$ be the quotient of $\Kb^{\tilde G_\l}(Z_\l)$
by its torsion $\Rb(\tilde G_\l)$-submodule.
The space $\Ub_\l$ is an associative algebra for the convolution
product $\star$ (with $Z_{12}=Z_{23}=Z_\l$).
It is proved in \cite{N3, Theorems 9.4.1 and 12.2.1} 
that the map $\xb^\pm_{ir}\mapsto x^\pm_{ir},$ 
$\kb^\pm_{ir}\mapsto k^\pm_{ir},$ 
extends uniquely to an algebra homomorphism
$\Phi_\l\,:\,\Ub\to\Ub_\l.$

\vskip3mm

\noindent{\bf Remark 6.3.} 
The morphism $\Phi_\lambda$ is not the one used by Nakajima, although the 
operators $h_{ir}$ in (6.2) and in \cite{N3, \S 9.2} 
are the same. The proof of Nakajima still works in our case : 
the relations \cite{N3, (1.2.8) and (1.2.10)} are checked in the appendix,
the relations involving only one vertex of the graph 
are proved as in \cite{N3, \S 11}, the Serre relations are proved as in 
\cite{N3, \S 10.4}.

\vskip3mm

\subhead 6.4\endsubhead
Recall that $\dim Q_{\a\l}=(\a|2\l-\a)$.
From the formula for $\Tc_{\a\l}$ in Section 4.5 we get 
$\Omega_{Q_{\a\l}}=q^{-d_{Q_{\a\l}}}$. Thus the hypothesis in Section 5.2 
are satisfied. Consider the anti-automorphism $\gamma_U=\phi_*D_{Z_\l}$ 
of $\Ub_\l$.

\vskip3mm

\noindent{\bf Lemma 6.5.} {\it We have $\Phi_\l\tau=\gamma_U\Phi_\l$.}

\vskip3mm

\noindent{\it Proof of 6.5.}
For any $\a'\in Q^+$ the Hecke correspondence
$C^{i\pm}_{\a'\l}$ is smooth and 
$$\Omega_{C^{i\pm}_{\a'\l}}=
q^{\mp f^i_{\a\l}}p^*\Omega_{Q_{\a\l}}\otimes p^*\Det(\Fc^{i*}_{\a\l})
\otimes (\pm\Lc)^{f^i_{\a\l}}_{|C^{i\pm}_{\a'\l}}$$
where $\a=\a'\mp\a_i$ (see the proof of 7.4).
Using the identities
$${p'}^*\Fc^{i\pm*}_{\a'\l}-p^*\Fc^{i\pm*}_{\a\l}=-q^{\mp 1}\Lc\and
\dim Z_{\a\a'\l}=\mp f^i_{\a\l}+d_{Q_{\a'\l}}+1,$$
and the commutation of the Serre-Grothendieck duality with closed embeddings
we get
$$\matrix
\gamma_U(x_{ir}^\pm)
&=\sum_{\a'}q^{-1}x^{i\mp}_{\a'\l}\star 
(\mp\Lc)_{|C^{i\mp}_{\a'\l}}^{\otimes f^{i\mp}_{\a'\l}-r}\hfill\cr\cr
&=\sum_{\a}
(\mp\Lc)_{|C^{i\mp}_{\a\l}}^{\otimes f^{i\mp}_{\a\l}-r}\star x^{i\mp}_{\a\l}
\hfill\cr\cr
&=x_{i,-r}^\mp.\hfill\cr
\endmatrix$$
\qed

\subhead 6.6\endsubhead
Put $\Wb_\l=\Kb^{\tilde G_\l}(L_\l)$,
$\Wb'_\l=\Kb^{\tilde G_\l}(Q_\l)$, and $\Rb_\l=\Rb(\tilde G_\l)$.
The $\Rb_\l$-modules $\Wb'_\l$, $\Wb_\l$ are free.
Thus $\Wb_\l$, resp. $\Wb'_\l$, may be viewed as $\Ub$-module via
the algebra homomorphism $\Ub\to\Ub_\l\to\End\,\Wb_\l$,
resp. $\Ub\to\End\,\Wb'_\l$, which is composed of $\Phi_\lambda$ 
and of the convolution product
$\star\,:\,\Ub_\lambda\boxtimes\Wb_\lambda\to\Wb_\lambda$,
resp. $\star\,:\,\Ub_\lambda\boxtimes\Wb'_\lambda\to\Wb_\lambda'$
(for $Z_{12}=Z_\l$ and $Z_{23}=Q_\l$, resp. $Z_{23}=L_\l$).
The varieties $L_{0\lambda}$ and $Q_{0\lambda}$ are
reduced to a point.
Let $[0]$ be their fundamental class in $\Kb$-theory. 
By \cite{N3, Propositions 12.3.2, 13.3.1} 
the $\Ub$-module $\Wb_\lambda$ 
is cyclic generated by $[0]$, and we have
$$\xb_i^+(z)\star [0]=0,\qquad\kb_i^\pm(z)\star[0]=q^{(\lambda\,|\,\a_i)}
\Wedge_{-1/z}\bigl((q^{-2}-1)\Wc_i\bigr)^\pm\otimes[0].$$

\subhead 6.7\endsubhead
Fix a semi-simple element $s=(t,\zeta)$ in $\tilde{G}_\l.$
Let $\la s\ra\subset \tilde{G}_\l$ be the Zariski closed subgroup generated
by $s$. Put 
$$\Wb_s=\Kb^{\la s\ra}(L_\l)_s,\quad
\Wb'_s=\Kb^{\la s\ra}(Q_\l)_s,\quad
\Ub_s=\Kb^{\la s\ra}(Z_\l)_s.$$
Let $\Phi_s\, :\,\Ub\to\Ub_s$
be the composition of $\Phi_\lambda$ and the specialization at $s$.
Consider the $\CC$-algebra
$$\Ub|_{q=\zeta}=\Ub\otimes_\AA\bigl(\AA/(q-\zeta)\bigr).$$
The spaces $\Wb_s$, $\Wb'_s$, are $\Ub|_{q=\zeta}$-modules. 
The $\Ub|_{q=\zeta}$-module $\Wb_s$
is called a standard module.
It is cyclic generated by $[0]$.

\vskip3mm

\head 7. The coproduct\endhead
\subhead 7.1\endsubhead
Fix $\l^1,\l^2\in P^+,$ such that $\l=\l^1+\l^2$. 
Put $G_{\l^1\l^2}=G_{\l^1}\times G_{\l^2}\times\CC^\times$,
$\Rb_{\l^1\l^2}=\Rb(G_{\l^1\l^2})$, and let 
$\bar\Rb_{\l^1\l^2}$ be the fraction field of $\Rb_{\l^1\l^2}$.
Set
$$\bar\Wb'_\l=\Wb'_\l\boxtimes_{\Rb_\l}\bar\Rb_{\l^1\l^2},\quad
\Wb'_{\l^1\l^2}=\Kb^{G_{\l^1\l^2}}(Q_{\l^1}\times Q_{\l^2}),\quad
\bar\Wb'_{\l^1\l^2}=\Wb'_{\l^1\l^2}\boxtimes_{\Rb_{\l^1\l^2}}\bar\Rb_{\l^1\l^2}.$$
Similarly we define $\bar\Wb_\l,\Wb_{\l^1\l^2},\bar\Wb_{\l^1\l^2}$
(using $L_\l$ instead of $Q_\l$), and
$\bar\Ub_\l,\bar\Ub_{\l^1\l^2}$ (using $Z_\l$ instead of $Q_\l$).
By \cite{N3, \S 7}, \cite{CG, \S 5.6} we have two Kunneth isomorphisms
$$\Wb'_{\lambda^1}\boxtimes\Wb'_{\lambda^2}\simto\Wb'_{\l^1\l^2},\quad
\Wb_{\lambda^1}\boxtimes\Wb_{\lambda^2}\simto\Wb_{\l^1\l^2}.$$
Let denote them by $\theta$. Set
$\tilde\Ub_{\l^1\l^2}=(\Ub_{\l^1}\boxtimes\Ub_{\l^2})
\boxtimes_{\Rb_{\l^1\l^2}}\bar\Rb_{\l^1\l^2}.$
We do not know if there is a Kunneth isomorphism
$$\Kb^{\tilde G_{\l^1}}(Z_{\l^1})\boxtimes
\Kb^{\tilde G_{\l^2}}(Z_{\l^2})\simeq
\Kb^{G_{\l^1\l^2}}(Z_{\l^1}\times Z_{\l^2}).$$
However, it is easy to see that the map
$\theta\,:\,\tilde\Ub_{\l^1\l^2}\to\bar\Ub_{\l^1\l^2}$ 
induced by the external tensor product is invertible : 
if $H=H_{\l^1}\times H_{\l^2}\times\CC^\times\subseteq G_{\l^1\l^2}$ 
is a maximal torus, then
$$\matrix
\bar\Rb(H)\boxtimes_{\bar\Rb_{\l^1\l^2}}\tilde\Ub_{\l^1\l^2}
&\simeq\bar\Kb^{H_{\l^1}\times\CC^\times}(Z_{\l^1})\boxtimes
\bar\Kb^{H_{\l^2}\times\CC^\times}(Z_{\l^2})\hfill\cr
&\simeq\bar\Kb^{H_{\l^1}\times\CC^\times}(Q_{\l^1}\times Q_{\l^1})\boxtimes
\bar\Kb^{H_{\l^2}\times\CC^\times}(Q_{\l^2}\times Q_{\l^2})\hfill\cr
&\simeq\bar\Kb^{H}(Q^2_{\l^1}\times Q^2_{\l^2})\hfill\cr
&\simeq\bar\Kb^{H}(Z_{\l^1}\times Z_{\l^2})\hfill\cr
&\simeq\bar\Rb(H)\boxtimes_{\bar\Rb_{\l^1\l^2}}\bar\Ub_{\l^1\l^2}.\hfill
\endmatrix$$
Here we have used the localization theorem, the identity
$$(Z_{\l^b})^{H_{\l^b}\times\CC^\times}=
(Q_{\l^b}\times Q_{\l^b})^{H_{\l^b}\times\CC^\times}\quad\roman{for}\quad
b=1,2,$$ 
the Kunneth formula for $Q_{\l^1}\times Q_{\l^2}$,
and \cite{CG, Theorem 6.1.22} 
(which is valid, althougth $G_{\l^1\l^2}$ is not simply connected).
Taking the invariants under the Weyl group we get the required invertibility.
This invertibility is not needed in the sequel.
Set also
$$\matrix 
\Tc_+=qE_+(\Vc^1,\Vc^2)+L(\Vc^1,q\Wc^2-\Vc^2)+
q^{-1}E_+(\Vc^2,\Vc^1)+q^{-2}L(\Vc^2,q\Wc^1-\Vc^1)\hfill\cr
\Tc'_+=qE(\Vc^1,\Vc^2)+L(\Vc^1,q\Wc^2-\Vc^2)+q^2L(q\Wc^1-\Vc^1,\Vc^2)
\hfill
\endmatrix$$
in $\Wb'_{\lambda^1\lambda^2},$
where $E_\pm(\Vc^1,\Vc^2)=\bigoplus_{i,j}n^\pm_{ij}{\Vc_i^1}^*\otimes\Vc_j^2$.

\proclaim{\bf Lemma 7.2} 
The class of the normal bundle of $Q_{\l^1}\times Q_{\l^2}$ in $Q_\l$ 
is $\Tc'_++q^2{\Tc'_+}^*$. In particular, the class
$\Omega'=\Wedge_{-1}({\Tc_+'}^*+q^{-2}\Tc_+')$ is well-defined and 
is invertible in $\bar\Wb'_{\lambda^1\lambda^2}.$
\endproclaim

\noindent{\it Proof of 7.2.}
Follows from Lemma 4.4 and Sections 4.5, 5.3.
\qed

\vskip3mm

\noindent
Let $\Delta'_{W'}\,:\,\bar\Wb'_\l\to\bar\Wb'_{\l^1\l^2}$
be the map induced by the pull-back 
$\kappa^*\,:\,\Wb'_\l\to\Wb'_{\l^1\l^2}$
(which is well defined, since $Q_{\l^1}\times Q_{\l^2}$ 
and $Q_\l$ are smooth).
Since $Z_\l$ is a closed subvariety of $(Q_\l)^2$,
the restriction with support with respect to the embedding 
$(Q_{\l^1})^2\times (Q_{\l^2})^2\hookrightarrow(Q_\l)^2$
gives a map
$$(\kappa\times\kappa)^*\,:\,\Kb^{G_\lambda\times\CC^\times}(Z_\lambda)\to
\Kb^{G_{\l^1}\times G_{\l^2}\times\CC^\times}
(Z_{\l^1}\times Z_{\l^2}).$$
Put $\Delta'_U=(1\boxtimes{\Omega'}^{-1})\star(\kappa\times\kappa)^*\,:\,
\bar\Ub_\l\to\bar\Ub_{\l^1\l^2},$
where $1\boxtimes{\Omega'}^{-1}$ is the pull-back of ${\Omega'}^{-1}$ 
by the second projection
$Z_{\l^1}\times Z_{\l^2}\to Q_{\l^1}\times Q_{\l^2}$. 
There is a unique linear map 
$Q\rightarrow P,$ $\alpha\mapsto \alpha^+$, such that
$f^{i+}_{\a\l}=(\a_i|\l-\a^+).$ 
By (5.4) and Lemma 7.2,
the class $\Wedge_{-1}(\Tc_+|_{Q_{\a^1\l^1}\times Q_{\a^2\l^2}})$
is well-defined and is invertible in the ring
$\bar\Kb^{G_{\l^1}\times G_{\l^2}\times\CC^\times}
(Q_{\a^1\l^1}\times Q_{\a^2\l^2})$.
Set
$$\Omega=\sum_{\a^1,\a^2}q^{(\a^1|\a^{2+}-\l^2)}
\Wedge_{-1}(\Tc_+|_{Q_{\a^1\l^1}\times Q_{\a^2\l^2}})
\in\bar\Wb'_{\lambda^1\lambda^2}.$$
The class $\Omega$ is invertible.
We put $\Delta^\diamond_{W'}=\Omega^{-1}\otimes\Delta'_{W'}$
(here $\otimes$ is the tensor product on $Q_{\l^1}\times Q_{\l^2}$).
Let $\delta_*\Omega^{\pm 1}$ be the image of $\Omega^{\pm 1}$ 
in $\bar\Ub_{\l^1\l^2}$.
We put $\Delta^\diamond_U=\delta_*\Omega^{-1}\star\Delta'_U\star\delta_*\Omega.$
Hereafter $\delta_*$ may be omitted.

\subhead 7.3\endsubhead
Let $\tilde R\in\tilde\Ub_{\l^1\l^2}$ be the element defined in Lemma 8.1.(iv),
and set $\bar R=\theta(\tilde R).$ 
We put $\Delta_{W'}=\bar R^{-1}\star\Delta^\diamond_{W'}$ 
and $\Delta_U=\bar R^{-1}\star\Delta^\diamond_U\star\bar R$.
Recall that we have anti-involutions $\gamma_U$ of $\bar\Ub_\l$ 
and $\bar\Ub_{\l^1\l^2}$ (see $\S 6.4$).
We set $\Delta^\gamma_U=\gamma_U\Delta_U\gamma_U$.

\proclaim{Theorem 7.4} The map 
$\Delta_U\,:\,\bar\Ub_\lambda\to\bar\Ub_{\lambda^1\lambda^2}$ satisfies
$$\Delta_U\Phi_\l=\theta(\Phi_{\l^1}\boxtimes\Phi_{\l^2})\Delta^\circ\and
\Delta^\gamma_U\Phi_\l=\theta(\Phi_{\l^1}\boxtimes\Phi_{\l^2})\Delta^\bullet.$$
\endproclaim

\proclaim{Lemma 7.5} Assume that $M,M',$ are smooth quasi-projective
$G$-varieties. Let $p\,:\,M\times M'\to M$ be the projection.
Fix a semisimple element $s\in G$ and a smooth closed $G$-subvariety
$X\subset M\times M'$. 
Put $\Nc=\Tc X-(p^*\Tc M)|_X$, and $\Nc^s=\Tc X^s-(p^*\Tc M^s)|_{X^s}$.
\hfill\break
(i) The element $\Wedge_{-1}(-\Nc^*|_{X^s}+\Nc^{s*})\in\Kb^{\la s\ra}(X^s)$ 
is well-defined.
Its image in $\Kb(X^s)$  under the evaluation map is still denoted by
$\Wedge_{-1}(-\Nc^*|_{X^s}+\Nc^{s*}).$ 
\hfill\break
(ii) For any $G$-bundle $\Ec$ on $X$, the bivariant localization morphism
$\rb\,:\,\Kb^{\la s\ra}(X)_s\to\Kb(X^s)$ 
defined in \cite{CG, \S 5.11} maps $\Ec$
to $\Ec|_{X^s}\otimes\Wedge_{-1}(-\Nc^*|_{X^s}+\Nc^{s*})$ 
(here $\otimes$ is the tensor product on $X^s$).
\endproclaim

\noindent{\bf Remark 7.6.} 
Fix $s=(t,\zeta)\in\tilde G_\lambda,$ $s^1=(t^1,\zeta)\in\tilde G_{\lambda^1},$ 
$s^2=(t^2,\zeta)\in\tilde G_{\lambda^2},$ such that $t=t^1\oplus t^2$ and
$t$ is semi-simple. Put
$\Ub_{s^1s^2}=\Kb^{\la s\ra}(Z_{\lambda^1}\times Z_{\lambda^2})_s,$
and let
$$\rb_s\,:\,\Ub_s\to\Kb(Z_\lambda^s),\quad
\rb_{s^1s^2}\,:\,\Ub_{s^1s^2}\to
\Kb(Z^{s^1}_{\lambda^1}\times Z^{s^2}_{\lambda^2})$$
be the bivariant localization maps.
These maps are invertible and commutes to the convolution product $\star$.
If $(\zeta^\ZZ\spec t^1)\cap\spec t^2=\emptyset$, then 
$Q_\l^s\simeq Q^{s^1}_{\l^1}\times Q^{s^2}_{\l^2}$, 
$Z_\l^s\simeq Z^{s^1}_{\l^1}\times Z^{s^2}_{\l^2}$, 
and the specialization
of $\Delta'_U$ at $s$ is well-defined and it coincides with
the map $\rb_{s^1s^2}^{-1}\cdot\rb_s$. 

\vskip3mm

\noindent{\it Proof of 7.5.} 
Claim $(i)$ is well-known, see \cite{CG, Proposition 5.10.3} for instance.
Claim $(ii)$ is immediate from the Koszul resolution
of $\Oc_X$ by sheaves of locally free $\Oc_{M\times M'}$-modules in a 
neighborhood, in $M\times M'$, of each point of $X^s$.\qed

\vskip3mm

\noindent{\it Proof of 7.4.} 
Assume that $s$ is generic. Then 
$\kappa(Q_{\lambda^1}^{s^1}\times Q_{\lambda^2}^{s^2})=Q_\lambda^s$ 
and $\Delta'_U$ specializes to the map $\rb_{s^1s^2}^{-1}\cdot\rb_s$ 
by Remark 7.6. 
Assume that $\a'=\a\pm\a_i\in Q^+$. Observe that 
$$\matrix
\Tc C^{i+}_{\a'\lambda}-\Tc Q_{\a\lambda}\boxtimes 1
&=1\boxtimes\Tc_{\alpha'\lambda}-\Tc^{i+}_{\alpha'\lambda}\hfill\cr
&=E(\Lc,\Vc')+L(\Lc,\Wc)-(q+q^{-1})L(\Lc,\Vc')+q^2\hfill\cr
&=q^2+q\Lc^*\otimes{p'}^*\Fc^i_{\alpha'\lambda},\hfill\cr\cr
\Tc C^{i-}_{\a'\lambda}-\Tc Q_{\a\lambda}\boxtimes 1
&=q^2+q\Lc\otimes{p'}^*\Fc^{i*}_{\alpha'\lambda}.\hfill
\endmatrix$$
We have
$$C^{i\pm}_{\a'\l}\cap(Z_{\l^1}\times Z_{\l^2})=
\bigsqcup_{\a=\a^1+\a^2}\bigl(C^{i\pm}_{{\a^1}'\l^1}\times\delta Q_{\a^2\l^2}
\sqcup\delta Q_{\a^1\l^1}\times C^{i\pm}_{{\a^2}'\l^2}\bigr),$$
where ${\a^b}'=\a^b\pm\a_i$.
Fix $\a^1,\a^2\in Q^+$ such that $\a=\a^1+\a^2$.
Take $M=Q_{\a\lambda}$, $M'=Q_{\a'\lambda}$, and $X=C^{i\pm}_{\a'\lambda}$
in Lemma 7.5. 
Let $\theta\,:\,\Ub_{s^1}\boxtimes\Ub_{s^2}\to\Ub_{s^1s^2}$
be the obvious map.
The element $\rb_{s^1s^2}^{-1}\cdot\rb_s(x^+_{ir})\in\Ub_{\s^1\s^2}$ 
is the image by $\theta$ of
$$\matrix
(x^+_{ir}\boxtimes 1)\otimes(-1)^{f^{i+}_{\l^2}}
\Det\bigl(\Lc\boxtimes\delta_*\Fc^{i+*}_{\l^2}\bigr)
\otimes\Wedge_{-1}\bigl(-q^{-1}\Lc\boxtimes\delta_*\Fc^{i*}_{\l^2}\bigr)+
\hfill\cr
+\Wedge_{-1}\bigl(-q^{-1}\delta_*\Fc^{i*}_{\l^1}\boxtimes\Lc\bigr)\otimes
(-1)^{f^{i+}_{\l^1}}\Det\bigl(\delta_*\Fc^{i+*}_{\l^1}\boxtimes\Lc\bigr)
\otimes(1\boxtimes x^+_{ir}).
\endmatrix$$
In this formula, the element 
$\Lc\boxtimes\delta_*\Fc^{i+*}_{\l^2}$ is identified with its
restriction to
$C^{i+}_{{\a^1}'\l^1}\times\delta Q_{\a^2\l^2}\subseteq
(Q_{\l^1})^2\times\delta(Q_{\l^2})$.
Thus $\Det$ is the maximal exterior power of a virtual bundle
on $C^{i+}_{{\a^1}'\l^1}\times\delta Q_{\a^2\l^2}.$
Similarly $\Wedge_{-1}\bigl(-q^{-1}\Lc\boxtimes\delta_*\Fc^{i*}_{\l^2}\bigr)$
is the $\Wedge_{-1}$ of a virtual bundle on
on $C^{i+}_{{\a^1}'\l^1}\times\delta Q_{\a^2\l^2}.$
It is well-defined by \cite{CG, Proposition 5.10.3} since $s$ is generic.
Here, $\otimes$ is the tensor product of sheaves 
on $C^{i+}_{{\a^1}'\l^1}\times\delta Q_{\a^2\l^2}$ (which is smooth).
In the same way, $\rb_{s^1s^2}^{-1}\cdot\rb_s(x^-_{ir})$ 
is the image by $\theta$ of
$$\matrix
(x^-_{ir}\boxtimes 1)\otimes 
(-1)^{f^{i-}_{\l^2}}\Det\bigl(-\Lc\boxtimes\delta_*\Fc^{i-*}_{\l^2}\bigr)
\otimes\Wedge_{-1}\bigl(-q^{-1}\Lc^*\boxtimes\delta_*\Fc^i_{\l^2}\bigr)+
\hfill\cr
+\Wedge_{-1}\bigl(-q^{-1}\delta_*\Fc^i_{\l^1}\boxtimes\Lc^*\bigr)\otimes
(-1)^{f^{i-}_{\l^1}}\Det\bigl(-\delta_*\Fc^{i-*}_{\l^1}\boxtimes\Lc\bigr)
\otimes(1\boxtimes x^-_{ir}).\hfill
\endmatrix$$
Since $s$ is generic, the class $\Omega\in\bar\Wb'_{\l^1\l^2}$
specializes to a class in $\Kb^{\la s\ra}(Q_{\l^1}\times Q_{\l^2})_s$.
Let $\Omega_{s^1s^2}$ be this class.
A direct computation gives the following identities in $\Ub_{s^1s^2}$
$$\matrix
\Omega^{-1}_{s^1s^2}\otimes(x_{ir}^\pm\boxtimes 1)\otimes\Omega_{s^1s^2}=
q^{\mp f^{i+}_{\l^2}}\Wedge_{-1}\bigl(q\Lc^*\boxtimes\Fc^{i+}_{\l^2}+
q^{-1}\Lc\boxtimes{\Fc^{i-}_{\l^2}}^*\bigr)
\otimes(x^\pm_{ir}\boxtimes 1),\hfill\cr\cr
\Omega^{-1}_{s^1s^2}\otimes(1\boxtimes x^\pm_{ir})\otimes\Omega_{s^1s^2}=
q^{\mp f^{i-}_{\l^1}}\Wedge_{-1}\bigl(q^{-1}\Fc^{i+}_{\l^1}\boxtimes\Lc^*+
q{\Fc^{i-}_{\l^1}}^*\boxtimes\Lc\bigr)
\otimes(1\boxtimes x^\pm_{ir}).\hfill
\endmatrix$$
Using (5.4) we get
$$\matrix
\Delta^\diamond_U(x^+_{ir})=\theta\Bigl(x^+_{ir}\boxtimes 1+q^{f^i_{\l^1}}
\Wedge_{-1}\bigl((q^{-1}-q)\Fc^i_{\l^1}\boxtimes\Lc^*\bigr)
\otimes(1\boxtimes x^+_{ir})\Bigr)\hfill\cr\cr
\Delta^\diamond_U(x^-_{ir})=
\theta\Bigl((x^-_{ir}\boxtimes 1)\otimes q^{f^i_{\l^2}}
\Wedge_{-1}\bigl((q^{-1}-q)(-\Lc)^*\boxtimes\Fc^i_{\l^2}\bigr)
+1\boxtimes x^-_{ir}\Bigr).\hfill
\endmatrix\leqno(7.7)$$
Using (7.7) it is proved in Lemma 8.1.(iii) that
$$\bar R^{-1}_{s^1s^2}\star(\Delta^\diamond_U\Phi_s)\star\bar R_{s^1s^2}=
\theta(\Phi_{s^1}\boxtimes\Phi_{s^2})\Delta^\circ.$$
We are done. The second identity follows from 
the first one and Lemma 6.5
since $\Delta^{\circ,\tau}=\Delta^\bullet$.
\qed

\vskip3mm

\subhead 7.8\endsubhead
Let $H\subset G_{\l^1}\times G_{\l^2}\times\CC^\times$ 
be the maximal torus of diagonal matrices.
Assume that the elements $s, s^1, s^2$ in $\S 7.3$ belong to $H.$ 
Put 
$$\Wb'_{\l,H}=\Kb^H(Q_\l),\quad
\Wb'_{\l^1\l^2,H}=\Kb^H(Q_{\l^1}\times Q_{\l^2}),$$
and idem for $\Wb_{\l,H}$, $\Wb_{\lambda^1\lambda^2,H}$.
The corresponding $\bar\Rb(H)$-vector spaces are overlined 
(i.e. we set $\bar\Wb'_{\l,H}=\bar\Kb^H(Q_\l)$, etc).
Let $\theta$ denote the Kunneth isomorphisms
$$\Wb_{\l^1,H}\boxtimes_{\Rb(H)}\Wb_{\l^2,H}\simto\Wb_{\l^1\l^2,H},\quad
\Wb'_{\l^1,H}\boxtimes_{\Rb(H)}\Wb'_{\l^2,H}\simto\Wb'_{\l^1\l^2,H}.$$
Consider the bilinear pairing
$$(\ |\ )\,:\,\Wb'_{\l,H}\boxtimes_{\Rb(H)}\Wb_{\l,H}\to\Rb(H),\,
\Ec\boxtimes\Fc\mapsto q_*(\Ec\otimes\Fc),$$
where $q$ is the projection to a point and 
$\otimes$ is the tensor product of sheaves on $Q^s_\l$. 
The following lemma is proved as in 
\cite{N3, Proposition 12.3.2 and Theorem 7.3.5}.

\proclaim{Lemma 7.9} (i) The $\Ub$-module 
$\Wb_{\l,H}$ is generated by $\Rb(H)\otimes [0]$.
\hfill\break
(ii) The $\Rb(H)$-modules $\Wb'_{\l,H}, \Wb_{\l,H}$ are free and
the pairing $(\, |\,)$ is perfect.
\endproclaim

\noindent
Let $(\ |\ )$ denote also the pairing with the scalars extended 
to the field $\bar\Rb(H)$, and the pairing between
$\bar\Wb'_{\l^1\l^2,H}$ and $\bar\Wb_{\l^1\l^2,H}.$
The automorphism 
$D_{L_\l}\boxtimes id_{\bar\Rb(H)}\,:\,\bar\Wb_{\l,H}\to\bar\Wb_{\l,H}$
is still denoted by $D_{L_\l}$. By $\S 5.2$ we have 
$$D_{Z_\l}(u)\star D_{L_\l}(m)= D_{L_\l}(u\star m),\qquad\forall u\in\Ub_\l,
\forall m\in\bar\Wb_{\l,H}.$$
Given a $\Ub_\l$-module $M$, let $M^\flat$ be its contragredient module,
that is $M^\flat=M^*$ as a vector space and
$(uf)(m)=f(\phi_*(u)m)$ for all $f\in M^\flat,$ $m\in M$, $u\in\Ub_\l$.
The symbols $\circ$ and $\bullet$ denote the tensor product of 
$\Ub$-modules relative to the coproduct $\Delta^\circ$ and $\Delta^\bullet$ 
respectively. To simplify let $\Delta_{W'}$ denote also the map
$\Delta_{W'}\boxtimes id_{\bar\Rb(H)}\,:\,\bar\Wb'_{\l,H}\to
\bar\Wb'_{\l^1\l^2,H}.$

\proclaim{Proposition 7.10} 
(i) The map $\Delta_{W'}\,:\,\bar\Wb'_{\l,H}\to\bar\Wb'_{\l^1\l^2,H}$ 
is invertible.\hfill\break
(ii) The map $\Delta_U$ is an algebra homomorphism and
$$\Delta_U(u)\star\Delta_{W'}(m')=\Delta_{W'}(u\star m'),
\quad\forall u,m'.$$
(iii) The pairing identifies $\Wb'_{\lambda,H}$ with the contragredient
module $\Wb^\flat_{\l,H}$.\hfill\break
(iv) Set $\Delta_W=D_{L_\l}({}^t\Delta_{W'}^{-1})D_{L_\l}$, where the transpose 
is relative to the pairing $(\ |\ )$. We have
$$\Delta^\gamma_U(u)\star\Delta_W(m)=\Delta_W(u\star m),
\quad\forall u,m.$$
(v) The map $\Delta_W$ is an embedding of $\Ub$-modules
$\Wb_{\l,H}\hookrightarrow\Wb_{\l^1\l^2,H}.$
\endproclaim

\vskip3mm

\noindent{\it Proof of 7.10.}
Claim (i) follows from the localization theorem since the fixpoint sets
$Q^H_\l$ and $Q^H_{\l^1}\times Q^H_{\l^2}$ are equal.
It suffices to check Claim (ii) on a dense subset of $\spec\Rb_{\l^1\l^2}$.
If $s$ is generic then $\Delta'_U=\rb_{s^1s^2}^{-1}\rb_s$ (see Remark 7.6).
Thus $\Delta_U$ is an algebra homomorphism, 
since the map $\rb_s$ commutes to the convolution product.
The case of $\Delta_{W'}$ is similar.
Claim (iii) means that
for any $m\in\Wb_{\l,H}$, $m'\in\Wb'_{\l,H}$, $u\in\Ub_\l$, we have 
$$(m'|u\star m)=(\phi_*(u)\star m'|m).$$ 
It suffices to check the two identities below :
$$(m'|u\star m)=(m'\star u|m)\and m'\star u=\phi_*(u)\star m'.$$
These identities are standard. The first one is the associativity 
of the convolution product, the second one is essentially the 
fact that $\phi_*$ is an anti-homomorphism.
Claim (iv) is a direct computation : 
for any $u,m,m'$ as above, we have 
$$\matrix
\bigl(\Delta_{W'}(m')|D_{L_\l}\Delta_W(u\star m)\bigr)
&=\bigl(m'|D_{L_\l}(u\star m)\bigr)\hfill\cr
&=\bigl(m'|D_{Z_\l}(u)\star D_{L_\l}(m)\bigr)\hfill\cr
&=\bigl(\gamma_U(u)\star m'|D_{L_\l}(m)\bigr)\hfill\cr
&=\bigl(\Delta_{W'}(\gamma_U(u)\star m')|D_{L_\l}\Delta_W(m)\bigr)\hfill\cr
&=\bigl(\Delta_U\gamma_U(u)\star\Delta_{W'}(m')|D_{L_\l}\Delta_W(m)\bigr)\hfill\cr
&=\bigl(\Delta_{W'}(m')|D_{Z_\l}\Delta_U^\gamma(u)\star D_{L_\l}\Delta_W(m)\bigr)\hfill\cr
&=\bigl(\Delta_{W'}(m')|D_{L_\l}(\Delta_U^\gamma(u)\star\Delta_W(m))\bigr).\hfill
\endmatrix$$
Since $\Wb_{\l,H}$ is a free $\Rb(H)$-module
the restriction of $\Delta_W$ to $\Wb_{\l,H}$ is injective. 
The $\Ub$-module $\Wb_{\l,H}$ is generated
by $\Rb(H)\otimes [0]$ and $\Delta_W([0])=[0]\boxtimes[0]$. Thus 
$\Delta_W(\Wb_{\l,H})\subseteq\Wb_{\l^1\l^2,H}.$
\qed

\subhead 7.11\endsubhead
We can now state the main result of this paper.
Using Theorem 7.4, Proposition 7.10, and the Kunneth isomorphisms
$$\theta^{-1}\,:\,
\Wb_{\l^1\l^2,H}\,\simto\,\Wb_{\l^1,H}\boxtimes_{\Rb(H)}\Wb_{\l^2,H},\quad
\theta^{-1}\,:\,
\Wb'_{\l^1\l^2,H}\,\simto\,\Wb'_{\l^1,H}\boxtimes_{\Rb(H)}\Wb'_{\l^2,H},$$
we get morphisms of $\Ub|_{q=\zeta}$-modules 
$$\theta^{-1}\Delta_W\,:\,\Wb_s\to\Wb_{s^1}\bullet\Wb_{s^2}\and
\theta^{-1}\Delta_{W'}\,:\,\Wb'_s\to\Wb'_{s^1}\circ\Wb'_{s^2}.$$

\proclaim{Theorem 7.12} 
Assume that $(\zeta^{-1-\NN}\spec t^1)\cap\spec t^2=\emptyset$.
Then the maps 
$$\theta^{-1}\Delta_W\,:\,\Wb_s\to\Wb_{s^1}\bullet\Wb_{s^2}\and
\theta^{-1}\Delta_{W'}\,:\,\Wb'_s\to\Wb'_{s^1}\circ\Wb'_{s^2}$$
are isomorphisms of $\Ub|_{q=\zeta}$-modules.
\endproclaim

\noindent{\bf Remark 7.13.} 
Let first recall the following standard facts
(see \cite{N3, Section 7.1} for instance).
Let $X$ be a smooth quasi-projective $G$-variety ($G$ a linear group)
with a finite partition into $G$-stable locally closed subsets
$X_i$, $i\in I$. Fix an order on $I$ such that the subset 
$\bigcup_{i'\leq i}X_{i'}\subset X$ is closed for all $i$.
Let $\kappa\,:\,Y\hookrightarrow X$ be the embedding of a smooth
closed $G$-stable subvariety such that
the intersections $Y_i=X_i\cap Y$, 
are the connected components of $Y$ (in particular, $Y_i$ is smooth). 
Assume that there is a $G$-invariant vector 
bundle map $\pi_i\,:\,X_i\to Y_i$ for each $i$.
Let $\Kb^G_{r,\top}$ be the complexified equivariant topological $K$-group
of degree $r$ ($r=0,1$). Assume also that $\Kb^G_{1,\top}(Y)=\{0\}$, 
$\Kb^G_{0,\top}(Y)=\Kb^G(Y)$ and $\Kb^G(Y)$ is a free $\Rb(G)$-module.
Let consider the vector bundle 
$\Nc_i=(\Tc X|_{Y_i})/(\Tc X_i|_{Y_i})$ on $Y_i$. 
For each $i$ fix a basis $(\Ec_{ij}; j\in J_i)$ of the space $\Kb^G(Y_i)$. 
Fix also an element $\bar\Ec_{ij}\in\Kb^G(\bar X_i)$ whose 
restriction to $X_i$ is $\pi^*_i\Ec_{ij}$. 
Then $\Kb^G_{1,\top}(X)=\{0\}$, $\Kb^G_{0,\top}(X)=\Kb^G(X),$
and the $\bar\Ec_{ij}$'s form a basis of $\Kb^G(X)$.
Moreover $\kappa^*(\bar\Ec_{ij})=\Wedge_{-1}(\Nc_i^*)\otimes\Ec_{ij}$ modulo
$\Kb^G(\bigcup_{i'<i}Y_{i'})$ (here $\otimes$ is the tensor product on $Y$).

\vskip3mm

\noindent{\it Proof of 7.12.} 
Consider the co-character
$$\gamma\,:\,\CC^\times\to\tilde G_\lambda,
\,z\mapsto(z\,id_{|\lambda^1|}\oplus id_{|\lambda^2|},1).$$
Let $\la s,\gamma\ra$ be the Zariski closed subgroup generated by $s$ 
and $\gamma(\CC^\times)$. We have
$\kappa(Q_{\lambda^1}^{s^1}\times Q_{\lambda^2}^{s^2})=
Q_\lambda^{\la s,\gamma\ra}$ (see Remark 7.6).
We claim that $\gamma$ gives a Byalinicki-Birula partition 
of the variety $Q_\l^s$ such that each piece is a $H$-equivariant vector bundle
over a connected component of $Q_\l^{\la s,\gamma\ra}$.
Fix $I$-graded vector spaces $V, W^1, W^2$, of dimension $\a,\l^1,\l^2$.
Given a triple $(B,p,q)$ representing a point $x\in Q_\lambda$ let
$V^1$ be the largest $B$-stable subspace contained in $q^{-1}(W^1)$.
Assume that $x$ is fixed by $s$. Let $Q(\rho)\subset Q_\lambda^s$ be the 
connected component containing $x$.
For any $z\in\CC^\times$ let $V(z)$ and $W(z)$ be defined as in Lemma 4.4. 
Assume that $(\zeta^{-1-\NN}\spec t^1)\cap \spec t^2=\emptyset$.
Then $V(z)\subset V^1$ for any $z\in\zeta^{-\NN}\spec(t^1)$.
In particular $p(W^1)\subset V^1$.
Thus, the subspace $V^1\oplus W^1\subset V\oplus W$ is stable by $B, p, q$.
The restriction of $(B,p,q)$ to $V^1\oplus W^1$ is a stable triple.
The projection of $(B,p,q)$ to $V/V^1\oplus W/W^1$ is stable either by 
the maximality of $V^1$.
Let $x^1\in Q_{\lambda^1}^{s^1}$, $x^2\in Q_{\lambda^2}^{s^2}$, 
be the classes of those triples.
We have $$\kappa(x^1,x^2)=\lim_{z\to 0}\gamma(z)x$$
(set $g(z)=z\,id_{V^1}\oplus id_S$ where $S$ is a $I$-graded vector space
such that $V=S\oplus V^1$,
and write the triple $\gamma(z)g(z)(B,p,q)$, which represents $\gamma(z)x$,
in a basis adapted to the splitting $V=S\oplus V^1$. 
Then, do the limit $z\to 0$). The claim is proved.

Consider the piece 
$$Q^+_{\rho^1\rho^2}=\{x\in Q_\lambda^s\,|\,
\lim_{z\to 0}\gamma(z)x\in Q_{\rho^1\rho^2}\},$$
where $Q_{\rho^1\rho^2}$ is the connected component 
$\kappa\bigl(Q(\rho^1)\times Q(\rho^2)\bigr)\subseteq Q_\lambda^{\la s,\gamma\ra}$.
By definition of $Q^+_{\rho^1\rho^2}$, the one-parameter
subgroup $\gamma$ acts on
$\Tc Q^+_{\rho^1\rho^2}|_{Q_{\rho^1\rho^2}}$ with non-negative weights. 
By Lemma 7.2 the class of the normal bundle of $Q_{\l^1}\times Q_{\l^2}$ in
$Q_\l$ is the sum of $\Tc_+'$ and $q^2{\Tc_+'}^*$.
It is easy to see that $\gamma$ acts on $\Tc_+'$
with negative weights, and on $q^2{\Tc_+'}^*$ with positive weights.
We get the following equality of classes in $\Kb^H(Q_{\rho^1\rho^2})$
$$\Tc Q_\lambda^s|_{Q_{\rho^1\rho^2}}-
\Tc Q^+_{\rho^1\rho^2}|_{Q_{\rho^1\rho^2}}=
\bigl(\Tc_+'|_{Q_{\rho^1\rho^2}}\bigr)^s.\leqno(7.14)$$
Here we use the notation : if $\Ec$ is the class of
a virtual $H$-bundle on $Q_{\rho^1\rho^2}$, then $\Ec^s$
is the class of the $s$-invariant part of the virtual bundle. 
Recall that 
$$Q_\l^{\la s,\gamma\ra}=\kappa(Q^{s^1}_{\l^1}\times Q^{s^2}_{\l^2})=
\bigsqcup_{\rho^1,\rho^2}Q_{\rho^1\rho^2}.$$
We can apply Remark 7.13 to the following situation 
$$X=Q^s_\l,\quad Y=\kappa(Q_{\l^1}^{s^1}\times Q_{\l^2}^{s^2}),\quad
I=\{(\rho^1,\rho^2)\},\quad
\{X_i\}=\{Q^+_{\rho^1\rho^2}\},\quad
\{Y_i\}=\{Q_{\rho^1\rho^2}\}$$
(see also \cite{N3, Theorem 7.3.5}). 
We get particular bases $\Bc_\l$ of $\Kb^H(Q^s_\l)$,
and $\Bc_{\l^1\l^2}$ of $\Kb^H(Q^{s^1}_{\l^1}\times Q^{s^2}_{\l^2})$.
In these bases, the $\bar\Rb(H)$-linear map
$$\Wedge_{-1}\bigl(({\Tc'}^*_+|_{Q_\lambda^{\la s,\gamma\ra}})^s\bigr)^{-1}
\otimes_1\kappa^*_1\,:\,
\bar\Kb^H(Q_\l^s)\to\bar\Kb^H(Q^{s^1}_{\l^1}\times Q^{s^2}_{\l^2})$$
(where $\kappa_1$ is the restriction of $\kappa$ to
$Q^{s^1}_{\l^1}\times Q^{s^2}_{\l^2}$, and $\otimes_1$ is the tensor product
on $\kappa(Q^{s^1}_{\l^1}\times Q^{s^2}_{\l^2})$)
is triangular unipotent by Remark 7.13 and (7.14).
Since $E_-(\Vc^1,\Vc^2)^*=E_+(\Vc^2,\Vc^1)$, the elements 
$$\Wedge_{-1}\bigl((\Tc_+|_{Q_{\rho^1\rho^2}})^s\bigr),\quad
\Wedge_{-1}\bigl(({\Tc'_+}^*|_{Q_{\rho^1\rho^2}})^s\bigr)\in
\bar\Kb^H(Q_{\rho^1\rho^2})$$
coincide up to the product by the class in $\Kb^H(Q_{\rho^1\rho^2})$
of an invertible sheaf (see $\S 7.1$ and $(5.4)$). 
Thus, the product by 
$$\Wedge_{-1}\bigl((\Tc_+|_{Q_{\rho^1\rho^2}})^s\bigr)^{-1}\otimes_1
\Wedge_{-1}\bigl(({\Tc'}^*_+|_{Q_{\rho^1\rho^2}})^s\bigr)$$
belongs to $\GL\bigl(\Kb^H(Q_{\rho^1\rho^2})\bigr)\subset
\GL\bigl(\bar\Kb^H(Q_{\rho^1\rho^2})\bigr)$.
In particular, the determinant of the $\bar\Rb(H)$-linear map
$$\Wedge_{-1}\bigl((\Tc_+|_{Q_\lambda^{\la s,\gamma\ra}})^s\bigr)^{-1}
\otimes_1\kappa^*_1\,:\,
\bar\Kb^H(Q_\l^s)\to\bar\Kb^H(Q^{s^1}_{\l^1}\times Q^{s^2}_{\l^2}),$$
respectively to $\Bc_\l$, $\Bc_{\l^1\l^2}$, 
is an invertible element of $\Rb(H)$. Let
$$\iota_\l\,:\,Q^s_\l\to Q_\l,\quad
\iota_{\l^1\l^2}\,:\,
Q^{s^1}_{\l^1}\times Q^{s^2}_{\l^2}\to Q_{\l^1}\times Q_{\l^2},$$
be the closed embeddings.
The localization theorem gives isomorphisms of $\bar\Rb(H)$-modules 
$$\iota_{\l*}\,:\,\bar\Kb^H(Q^s_\l)\simto\bar\Wb'_{\l,H},\quad
\iota_{\l^1\l^2*}\,:\,\bar\Kb^H(Q^{s^1}_{\l^1}\times Q^{s^2}_{\l^2})
\simto\bar\Wb'_{\l^1\l^2,H}.$$ 
Moreover, $\iota_{\l*}\bigl(\Kb^H(Q^s_\l)\bigr)\subset\Wb'_{\l,H}$
are free $\Rb(H)$-submodules of $\bar\Wb'_{\l,H}$ of maximal rank such that
$\iota_{\l*}\bigl(\Kb^H(Q^s_\l)\bigr)_s=\Wb'_s$. 
Thus, the basis $\iota_{\l*}(\Bc_\l)$ differs from any basis of $\Wb'_{\l,H}$
by the action of an $\bar\Rb(H)$-linear operator whose determinant 
is regular and non-zero at $s$. Idem for 
$\iota_{\l^1\l^2*}(\Bc_{\l^1\l^2})\subset\Wb'_{\l^1\l^2,H}$.
Set 
$$\bar f=\Omega^{-1}\otimes\kappa^*\,:
\,\bar\Wb'_{\lambda,H}\to\bar\Wb'_{\lambda^1\lambda^2,H}.$$
By \S 7.1 we have
$$\iota_{\l^1\l^2*}^{-1}\circ\bar f\circ\iota_{\l*}
=\sum_{\rho^1,\rho^2}A_{\rho^1\rho^2}\otimes_1
\Wedge_{-1}\bigl(\Tc_+|_{Q_{\rho^1\rho^2}}\bigr)^{-1}
\otimes_1\kappa^*_1,$$
for some element $A_{\rho^1\rho^2}\in\bar\Kb^H(Q_{\rho^1\rho^2})$.
Moreover, the product by 
$$A_{\rho^1\rho^2}\otimes_1
\Wedge_{-1}\bigl((\Tc_+|_{Q_{\rho^1\rho^2}})^s\bigr)\otimes_1
\Wedge_{-1}\bigl(\Tc_+|_{Q_{\rho^1\rho^2}}\bigr)^{-1}$$
is an invertible operator in $\GL\bigl(\bar\Kb^H(Q_{\rho^1\rho^2})\bigr)$
whose determinant, in $\bar\Rb(H)$, is regular and non-zero at $s$.
The element $\bar R\in\bar\Ub_{\l^1\l^2}$ is unipotent by Lemma 8.1.(v).
Thus the determinant of the $\bar\Rb(H)$-linear map
$$\Delta_{W'}\,:\,\bar\Wb'_{\l,H}\to\bar\Wb'_{\l^1\l^2,H},$$
with respect to $\iota_{\l*}(\Bc_\l),$ 
$\iota_{\l^1\l^2*}(\Bc_{\l^1\l^2})$
is regular and non zero at $s$. 
By Proposition 7.10.(iv), (v) we have
$\Delta_{W'}(\Wb'_{\l,H})\subset\Wb'_{\l^1\l^2,H}$.
Thus, the map 
$\theta^{-1}\Delta_{W'}\,:\,\Wb'_{\l,H}\to
\Wb'_{\l^1,H}\boxtimes_{\Rb(H)}\Wb'_{\l^2,H}$
specializes to an isomorphism 
$\Wb'_s\to\Wb'_{s^1}\circ\Wb'_{s^2}$ 
of $\Ub|_{q=\zeta}$-modules.
The other claim is due to the following easy fact. 
Consider the tensor product of $\Ub_\l$-modules
$M_1\boxtimes M_2$ relative to the coproduct $\Delta_U^\gamma$. If the map
$\Delta_M\,:\,M\to M_1\boxtimes M_2$ is an isomorphism of $\Ub_\l$-modules,
then the map 
${}^t\Delta_M^{-1}\,:\,M^\flat\to M_1^\flat\boxtimes_\phi M_2^\flat$ 
is an isomorphism of $\Ub_\l$-modules either, where
$M_1^\flat\boxtimes_\phi M_2^\flat$ is the tensor product relative
to the coproduct $\phi_*\Delta^{\gamma}_U\phi_*$.
\qed

\vskip3mm

\subhead 7.15\endsubhead
For any $\a\in\CC^\times$ and any $k\in I$, let $V_\zeta(\omega_k)_\a$ 
be the simple
finite dimensional $\Ub_{|q=\zeta}$-module with the $j$-th 
Drinfeld polynomial $(z-\zeta^{-1}\a)^{\delta_{jk}}$. 
By \cite{N3, Theorem 14.1.2} we have $V_\zeta(\omega_k)_\a=\Wb_s$ if $\lambda=\omega_k$
and $s=(\a,\zeta)\in\tilde G_\lambda.$
Theorem 7.12 has the following corollaries 
which were conjectured in \cite{N3} (in a less precise form) and in \cite{AK}
(for all types) respectively.

\vskip3mm

\noindent{\bf Corollary 7.16.} {\it 
The standard modules are the tensor products of the modules
$$V_\zeta(\omega_{i_1})_{\a_j\zeta^{\tau_1}}\circ\cdots\circ 
V_\zeta(\omega_{i_n})_{\a_j\zeta^{\tau_n}}$$
such that $\tau_1\geq\tau_2\geq...\geq\tau_n$, $j=1,2,...,r$, and the 
complex numbers $\a_j$ are distincts modulo $\zeta^\ZZ$.} 
\qed

\proclaim{Corollary 7.17} The $\Ub$-module
$V_\zeta(\omega_{i_1})_{\zeta^{\tau_1}}\circ\cdots\circ
V_\zeta(\omega_{i_n})_{\zeta^{\tau_n}}$
is cyclic if $\tau_1\geq\tau_2\geq...\geq\tau_n$ 
and is cocyclic if $\tau_1\leq\tau_2\leq...\leq\tau_n.$
Moreover it is generated (resp. cogenerated) by the tensor product of highest
weight vectors.
\endproclaim

\noindent{\bf Remark 7.18.} Corollary 7.16 is false if $\zeta$ is
a root of unity. 

\vskip3mm

\noindent{\bf Remark 7.19.} 
The module $\Wb_\lambda$ is presumably isomorphic to the Weyl module
introduced in \cite{K}, and studied in \cite{CP}. 
This statement is related to the flatness of the Weyl modules 
over the ring $\Rb_\lambda$.
This is essentially equivalent to the conjecture in \cite{CP}.

\vskip3mm

\head 8. The $R$-matrix\endhead
Fix $\l^1,\l^2\in P^+,$ and put $\l=\l^1+\l^2$. 
For any subset $T\subset\CC$ we put $|T|=\{|z|\,:\,z\in T\}$.
If $T, T'\subset\CC$ we write $|T|<|T'|$ if and only if
$t<t'$ for all $(t,t')\in |T|\times |T'|$.
Let $S$ be the set of pairs of semi-simple elements
$(s^1,s^2)\in\tilde G_{\l^1}\times\tilde G_{\l^2}$ such that
$s^1=(t^1,\zeta)$, $s^2=(t^2,\zeta)$ and
$|\spec\rho^1(s^1)^{-1}|<1<|\spec\rho^2(s^2)^{-1}|$ for all
$\rho^1,\rho^2$ such that $Q_{\rho^1\rho^2}$ is a non-empty connected
component of $Q_{\l^1}^{s^1}\times Q_{\l^2}^{s^2}$. 
The projection of $S$ to $\spec\Rb_{\l^1\l^2}$
is a Zariski-dense subset.
For any $s^1,s^2,$ we put $s=(t^1\oplus t^2,\zeta)$.
As usual, we put $\rho=\sum_{i\in I}\omega_i\in P^+$.

\proclaim{Lemma 8.1} Assume that $(s^1,s^2)\in S$. \hfill\break
(i) If $m^+=\xb^+_{i_1r_1}\cdots\xb^+_{i_kr_k},$
$m^-=\xb^-_{j_1s_1}\cdots\xb^-_{j_ks_k}\in\Ub,$
then $(\Phi_{s^1}\boxtimes\Phi_{s^2})(T_{n\rho}^{[2]}(m^-\boxtimes m^+))$ 
goes to zero when $n\to\infty$.
\hfill\break
(ii) The sequence $(\Phi_{s^1}\boxtimes\Phi_{s^2})(R_{2n\rho})$
admits a limit in $\Ub_{s^1}\boxtimes\Ub_{s^2}$ when $n$ goes to $\infty$.
This limit, denoted by $\tilde R_{s^1s^2}$, is an invertible element.
\hfill\break
(iii) Put $\bar R_{s^1s^2}=\theta(\tilde R_{s^1s^2})$. Then
$\bar R_{s^1s^2}\star\theta(\Phi_{s^1}\boxtimes\Phi_{s^2})
\Delta^\circ\star\bar R^{-1}_{s^1s^2}=\Delta^\diamond_U\Phi_s.$
\hfill\break
(iv) There is a unique invertible element $\tilde R\in\tilde\Ub_{\l^1\l^2}$
which specializes to $\tilde R_{s^1s^2}$ for any $(s^1,s^2)\in S.$ 
\hfill\break
(v) The element $\tilde R\in\tilde\Ub_{\l^1\l^2}$ is unipotent.
\endproclaim

\noindent{\it Proof of 8.1.}
In Part (i) we can assume that $k=1$, i.e. $m^-=x^-_{j,s}$ and $m^+=x^+_{i,r}$.
Let $\Lc$ be the virtual bundle on $Q_\l\times Q_\l$ introduced in \S 6.1. 
By (6.2) we have 
$$(\Phi_{\l^1}\boxtimes\Phi_{\l^2})(T_{n\rho}^{[2]}(m^-\boxtimes m^+))=
\pm\bigl(\Lc^{s+n}\boxtimes\Lc^{r-n}\bigr)\otimes 
\bigl(x^-_{j0}\boxtimes x^+_{i0}\bigr),$$
where $\otimes$ is the tensor product on $Q_\l\times Q_\l$.
For any $n>0$, let $\Lc^n_{s^b}$, $b=1,2$, be the image in $\Ub_{s^b}$ 
of the restriction of $\Lc^n$ to $Z_\l$.
We want to prove that 
$\lim_{n\to\infty}\Lc^n_{s^1}\boxtimes\Lc^{-n}_{s^2}=0.$
Fix $\rho^1,\rho^2$, such that the subset
$Q_{\rho^1\rho^2}\subset Q_{\l^1}^{s^1}\times Q_{\l^2}^{s^2}$ is non empty. 
The set of the eigenvalues of $s^b$ acting on
the bundle $\Vc|_{Q(\rho^b)}$ is the spectrum
of the semi-simple element $\rho^b(s^b)^{-1}$. 
For any $\a\in\spec\rho^b(s^b)^{-1}$, let
$\Vc_{\rho^b}^\a\subseteq\Vc|_{Q(\rho^b)}$
be the corresponding eigen-sub-bundle.
The image of $\Lc^n_{s^1}\boxtimes\Lc^{-n}_{s^2}$ 
by $\rb_{s^1}\boxtimes\rb_{s^2}$ is the restriction to
$Z_{\l^1}^{s^1}\times Z_{\l^2}^{s^2}$ of the product of
$$\sum_{\rho^1,\rho^2}\sum_{\a,\b}(\a/\b)^n
\bigl(\Vc^\a_{\rho^1}-{\Vc'}^\a_{\rho^1}\bigr)^{\otimes n}
\boxtimes\bigl(\Vc^{\b*}_{\rho^2}-{\Vc'}^{\b*}_{\rho^2}\bigr)^{\otimes n}$$
by a constant which does not depend on $n$.
Let us recall that if $\Ec$ is a line bundle on a smooth variety $X$, 
then the element $\Ec-1$ is a nilpotent element in the ring $\Kb(X)$
(see \cite{CG, Proposition 5.9.4} for instance).
Since $(s^1,s^2)\in S$, we are done.

Claim (ii) follows from Claim (i), from the formula 
$$R_{n\rho}=T_{(n-1)\rho}^{[2]}(R_\rho)T_{(n-2)\rho}^{[2]}(R_\rho)\cdots R_\rho,$$
and from \cite{D, Theorem 4.4(2)}, applied to the partial $R$-matrix $R_\rho$. 

We will prove Claim (iii) as in \cite{KT, Appendix B}. 
Fix $r\in\ZZ$. Since
$R_{2n\rho}\cdot\Delta^\circ(\xb^+_{ir})\cdot R_{2n\rho}^{-1}=
\Delta^\circ_{2n\rho}(\xb^+_{ir})$, see $\S 3$, the limit
$\ell im=\lim_{n\to\infty}(\Phi_{s^1}\boxtimes\Phi_{s^2})
\Delta^\circ_{2n\rho}(\xb_{ir}^+)$ is well-defined by Claim (ii). 
Let us prove that the series
$x^+_{ir}\boxtimes 1+\sum_{s\geq 0}k^+_{is}\boxtimes x^+_{i,r-s}$
converges, and coincides with this limit. We will assume that
$r\geq 0$, the case $r<0$ being very similar. The element 
$\Delta^\circ(\xb^+_{ir})-\xb^+_{ir}\boxtimes 1-\kb_i\boxtimes\xb^+_{ir}$
is a linear combination of monomials of the form
$$\Bigl(\vec\Prod_{u=1}^a\xb_{j_us_u}^-\Bigr)
\Bigl(\Prod_{u=1}^b\kb_{h_uq_u}^+\Bigr)\boxtimes
\Bigl(\vec\Prod_{u=a}^0\xb_{i_ur_u}^+\Bigr),$$
where 
$\vec\Prod_{u=1}^a$ is the ordered product 
(the term corresponding to $u=1$ is on the left),
$1\leq s_u\leq r\geq r_u\geq 0$,
$\sum_u\a_{i_u}-\sum_u\a_{j_u}=\a_i$,
and $\sum_us_u+\sum_ur_u+\sum_uq_u=r$,
see \cite{D, Theorem 4.4(3)} and its proof in
\cite {DD, \S 3.5}.
Since $\Delta^\circ_{2n\rho}(\xb^+_{ir})=
T^{[2]}_{2n\rho}\Delta^\circ(\xb^+_{i,r+2n}),$
the element 
$\Delta^\circ_{2n\rho}(\xb^+_{ir})-\xb^+_{ir}\boxtimes 1-\kb_i\boxtimes\xb^+_{ir}$
is, thus, a linear combination of monomials of the form 
$$\Bigl(\vec\Prod_{u=1}^a\xb_{j_u,s_u+2n}^-\Bigr)
\Bigl(\Prod_{u=1}^b\kb_{h_uq_u}^+\Bigr)\boxtimes
\Bigl(\vec\Prod_{u=a}^0\xb_{i_u,r_u-2n}^+\Bigr),\leqno(8.2)$$
where 
$0\leq s_u\leq r+2n\geq r_u\geq 0$,
$\sum_u\a_{i_u}-\sum_u\a_{j_u}=\a_i$,
and $\sum_us_u+\sum_ur_u+\sum_uq_u=r+2n$.
By Claim (i) the image of the monomials (8.2) by $\Phi_{s^1}\boxtimes\Phi_{s^2}$
cannot contribute to $\ell im$ unless $a=0$ 
(use the relation $\kb^+_{ir}=[\xb_{i0}^+,\xb_{ir}^-]$, and the inequalities
$\sum_u(s_u+2n)+\sum_uq_u\geq 2na$, $\sum_u(r_u-2n)\leq r-2na$).
Then, $i_0=i$. Let us consider the monomials
$$[(h_u),(q_u),r_0]=\Bigl(\Prod_u\kb_{h_uq_u}^+\Bigr)\boxtimes\xb_{ir_0}^+,
\quad\text{where}\quad r_0,q_u\geq 0
\quad\text{and}\quad r_0+\sum_uq_u=r.$$ 
We have $\Delta^\circ(\xb^+_{ir})=
A(\xb^+_{ir}\boxtimes 1-\kb_{i0}^+\boxtimes\xb^+_{ir})A^{-1},$
where $A=T_{-r\rho}^{[2]}(R_{r\rho})$.
By (3.1) the partial $R$-matrix $R_{r\rho}$
is a sum of monomials in the elements 
$\xb_{js}^-\boxtimes\xb^+_{j,-s}$ with $s\geq 1$.
Moreover, the coefficients of
$\xb_{i,-r_0}^-\boxtimes\xb^+_{ir_0}$ in $A$, $A^{-1}$
are respectively $c_1,\bar c_1$ 
(because the coefficients of
$\xb_{i,r-r_0}^-\boxtimes\xb^+_{i,r_0-r}$ in $R_{r\rho}$, $R_{r\rho}^{-1}$
are respectively $c_1,\bar c_1$, see \S 3).
Thus, the coefficient of $[(h_u),(q_u),r_0]$ in
$\Delta^\circ(\xb^+_{ir})-\xb^+_{ir}\boxtimes 1-\kb_{i0}^+\boxtimes\xb^+_{ir}$
and in 
$$(c_1\xb^-_{i,-r_0}\xb^+_{ir}+
\bar c_1\xb^+_{ir}\xb^-_{i,-r_0})\boxtimes\xb^+_{ir_0}=
(q-q^{-1})[\xb^+_{ir},\xb^-_{i,-r_0}]\boxtimes\xb^+_{ir_0}
=\kb^+_{i,r-r_0}\boxtimes\xb^+_{ir_0}$$
coincide. We are done.
A similar argument gives the equality (and the convergence of both sides)
$$\lim_{n\to\infty}(\Phi_{s^1}\boxtimes\Phi_{s^2})
\Delta^\circ_{2n\rho}(\xb_{ir}^-)=
1\boxtimes x^-_{ir}+\sum_{s\geq 0}x^-_{i,r+s}\boxtimes k^-_{i,-s},$$
for all $r\in\ZZ$.
By (6.2) the $r$-th Fourier coefficient in
$$x^+_i(z)\boxtimes 1+q^{f^i_{\l^1}}\Wedge_{-1}\bigl((q^{-1}-q)\,z^{-1}\,
\Fc^i_{\l^1}\bigr)^+\boxtimes x^+_i(z)$$
is $x^+_{ir}\boxtimes 1+\sum_{s\geq 0}k^+_{is}\boxtimes x^+_{i,r-s}$.
Similarly, $1\boxtimes x^-_{ir}+\sum_{s\geq 0}x^-_{i,r+s}\boxtimes k^-_{i,-s}$
is the $r$-th Fourier coefficient in
$$\bigl(x^-_i(z)\boxtimes 1\bigr)\otimes q^{f^i_{\l^2}}
\Wedge_{-1}\bigl((q^{-1}-q)z^{-1}\Fc^i_{\l^2}\bigr)^-+1\boxtimes x^-_i(z).$$
Since $(s^1,s^2)\in S$ the class
$\Wedge_{-1}\bigl((q^{-1}-q)\Fc^i_{\l^1}\boxtimes\Lc^*\bigr)$ is well-defined
in $\Ub_{s^1}\boxtimes\Ub_{s^2}$.
Recall that, by (7.7), the element $\Delta^\diamond_U(x^\pm_{ir})$ 
is the image by $\theta$ of the $r$-th Fourier coefficient in
$$\matrix
x^+_i(z)\boxtimes 1+
q^{f^i_{\l^1}}\Wedge_{-1}\bigl((q^{-1}-q)\Fc^i_{\l^1}\boxtimes\Lc^*\bigr)
\otimes\bigl(1\boxtimes x^+_i(z)\bigr),\hfill\cr\cr
\bigl(x^-_i(z)\boxtimes 1\bigr)\otimes q^{f^i_{\l^2}}
\Wedge_{-1}\bigl((q^{-1}-q)(-\Lc)^*\boxtimes\Fc^i_{\l^2}\bigr)
+1\boxtimes x^-_i(z).\hfill
\endmatrix$$
Since $\Lc\boxtimes x_i^\pm(z)=\pm zx_i^\pm(z)$ by (6.2), we get
$$\lim_{n\to\infty}\theta(\Phi_{s^1}\boxtimes\Phi_{s^2})
\Delta^\circ_{2n\rho}(\xb_{ir}^\pm)=
\Delta^\diamond_U(x^\pm_{ir}).\leqno(8.3)$$ 
We are done because the algebra $\Ub$ is generated 
by the elements $\xb_{ir}^\pm$.

Let us prove Part (iv).
The map $\Delta'_{W'}\,:\,\bar\Wb'_\l\to\bar\Wb'_{\l^1\l^2}$ 
is the bivariant localization map
(associated to the embedding
$(Q_{\l^1}\times Q_{\l^2})\times pt\hookrightarrow Q_\l\times pt$).
In particular, it is invertible.
The map $\Delta^\diamond_{W'}$ is invertible either. 
Thus $\bar\Wb'_{\l^1\l^2}$ can be viewed as a $\Ub$-module in two different 
ways : via $\Delta^\diamond_{W'}$ or via the coproduct $\Delta^\circ$ and  
the Kunneth isomorphism
$\bar\Wb'_{\l^1\l^2}\simeq
\bar\Wb'_{\l^1}\boxtimes_{\bar\Rb_{\l^1\l^2}}\bar\Wb'_{\l^2}$.
These representations of $\Ub$ are denoted by 
${}^\diamond\bar\Wb'_{\l^1\l^2}$ and ${}^\circ\bar\Wb'_{\l^1\l^2}$. 
Both $\Ub$-modules are finite dimensional $\bar\Rb_{\l^1\l^2}$-vector spaces.
They are simple by \cite{N3, Theorem 14.1.2}, 
since the tensor product of two generic simple 
finite-dimensional $\Ub$-modules is still simple. 
They have also the same Drinfeld polynomials (see \S 6.6).
Thus they are isomorphic. Obviously,
$$\Hom_{\bar\Rb_{\l^1\l^2}}
({}^\circ\bar\Wb'_{\l^1\l^2},{}^\diamond\bar\Wb'_{\l^1\l^2})=
\End_{\bar\Rb_{\l^1\l^2}}(\bar\Wb'_{\l^1\l^2}).$$
Let $H_{\l^b}\subset G_{\l^b}$ be a maximal torus, $b=1,2$.
Put $H=H_{\l^1}\times H_{\l^2}\times\CC^\times$. Then 
$$(Z_{\l^b})^{H_{\l^b}\times\CC^\times}=
(Q_{\l^b}\times Q_{\l^b})^{H_{\l^b}\times\CC^\times}.$$
The localization theorem and the Kunneth formula, give an isomorphism 
$$\bar\Kb^{H_{\l^b}\times\CC^\times}(Z_{\l^b})\simeq
\End_{\bar\Rb(H_{\l^b}\times\CC^\times)}
\bigl(\bar\Kb^{H_{\l^b}\times\CC^\times}(Q_{\l^b})\bigr).$$
Moreover, \cite{CG, Theorem 6.1.22} gives 
$$\matrix
\bar\Rb(H)\boxtimes_{\bar\Rb_{\l^1\l^2}}\bar\Wb'_{\l^1\l^2}\simeq
\bar\Kb^{H_{\l^1}\times\CC^\times}(Q_{\l^1})\boxtimes
\bar\Kb^{H_{\l^2}\times\CC^\times}(Q_{\l^2})\hfill\cr
\bar\Rb(H)\boxtimes_{\bar\Rb_{\l^1\l^2}}\tilde\Ub_{\l^1\l^2}\simeq
\bar\Kb^{H_{\l^1}\times\CC^\times}(Z_{\l^1})\boxtimes
\bar\Kb^{H_{\l^2}\times\CC^\times}(Z_{\l^2}).\hfill\cr
\endmatrix
$$
Thus,
$$\bar\Rb(H)\boxtimes_{\bar\Rb_{\l^1\l^2}}\Hom_{\bar\Rb_{\l^1\l^2}}
({}^\circ\bar\Wb'_{\l^1\l^2},{}^\diamond\bar\Wb'_{\l^1\l^2})\simeq
\bar\Rb(H)\boxtimes_{\bar\Rb_{\l^1\l^2}}\tilde\Ub_{\l^1\l^2}.$$
Taking the invariants by the Weyl group of $G_{\l^1}\times G_{\l^2}$, we get
an isomorphism
$$\Hom_{\bar\Rb_{\l^1\l^2}}
({}^\circ\bar\Wb'_{\l^1\l^2},{}^\diamond\bar\Wb'_{\l^1\l^2})\simeq
\tilde\Ub_{\l^1\l^2}.$$
Let $\tilde R\in\tilde\Ub_{\l^1\l^2}$ be the element corresponding to
the isomorphism of $\Ub$-modules
${}^\circ\bar\Wb'_{\l^1\l^2}\,\simto\,{}^\diamond\bar\Wb'_{\l^1\l^2}$
which is the identity on the highest weight vectors.
Claim (iv) follows from Claim (iii),
since $\tilde R_{s^1s^2}$ intertwines the specialized modules, whenever
it is defined and invertible.

Claim (v) is immediate, since $R_{2n\rho}$
is a sum of monomials in the elements 
$\xb_{ir}^-\boxtimes\xb^+_{i,-r}$, $i\in I$, $r\geq 1$.
\qed

\vskip3mm

\head A. Appendix\endhead
Let us check that the operators introduced in Section 6 still satisfy the 
Drinfeld relations. As indicated in Remark 6.3 it is sufficient to check
\cite{N3, (1.2.8) and (1.2.10)}. By \cite{N3, Section 10.2} the proof
of the first relation is reduced to the equality
$$\matrix
(q^{-1}\Vc^2_l/\Vc^3_l)^{f^{l-}_{\a^4\l}}\otimes
(q^{-1}\Vc^3_k/\Vc^4_k)^{f^{k+}_{\a^3\l}}\otimes
x^{l-}_{\a^4\l}\otimes x^{k+}_{\a^3\l}=\hfill\cr
=(q^{-1}\Vc^2_l/\Vc^3_l)^{f^{l-}_{\a^3\l}}\otimes
(q^{-1}\Vc^3_k/\Vc^4_k)^{f^{k+}_{\a^2\l}}\otimes
x^{k+}_{\a^2\l}\otimes x^{l-}_{\a^3\l},\hfill
\endmatrix\leqno(A.1)$$
where $\a^2,\a^3,\a^4\in Q^+,$ 
are such that $\a^4=\a^3-\a_k$,
$\a^2=\a^3+\a_l$, and $k\neq l$. 
Then $(A.1)$ follows from 
$$\matrix
x^{l-}_{\a^4\l}=(-1)^{n_{lk}^-}x^{l-}_{\a^3\l}
\otimes(q^{-1}\Vc_k^3/\Vc^4_k)^{n^-_{lk}},\qquad\hfill&
f^{l-}_{\a^4\l}=f^{l-}_{\a^3\l}-n^-_{lk},\hfill\cr
x^{k+}_{\a^2\l}=(-1)^{n_{kl}^+}x^{k+}_{\a^3\l}
\otimes(q^{-1}\Vc_l^2/\Vc^3_l)^{-n^+_{kl}},\hfill&
f^{k+}_{\a^2\l}=f^{k+}_{\a^3\l}+n^+_{kl},\hfill
\endmatrix$$
and the identity $n^+_{kl}=n^-_{lk}$.
By \cite{N3, Section 10.3} the proof of the second relation 
is reduced to the equality
$$\matrix
(q^{-1}\Vc^3_l/\Vc^2_l)^{a_{kl}}\otimes
(q^{-1}\Vc^3_k/\Vc^4_k)^{f^{k+}_{\a^2\l}}\otimes
(q^{-1}\Vc^3_l/\Vc^2_l)^{f^{l+}_{\a^3\l}}\otimes
x^{k+}_{\a^2\l}\otimes x^{l+}_{\a^3\l}=\hfill\cr
=(-1)^{a_{kl}}(q^{-1}\Vc^3_k/\Vc^4_k)^{a_{kl}}\otimes
(q^{-1}\Vc^3_l/\Vc^2_l)^{f^{l+}_{\a^4\l}}\otimes
(q^{-1}\Vc^3_k/\Vc^4_k)^{f^{k+}_{\a^3\l}}\otimes
x^{l+}_{\a^4\l}\otimes x^{k+}_{\a^3\l},\hfill
\endmatrix\leqno(A.2)$$
where $\a^2,\a^3,\a^4\in Q^+,$ are such that $\a^2=\a^3-\a_l$,
$\a^4=\a^3-\a_k$, and $k\neq l$. Then $(A.2)$ follows from 
$$\matrix
x^{l+}_{\a^4\l}=(-1)^{n_{lk}^+}
x^{l+}_{\a^3\l}\otimes(q^{-1}\Vc_k^3/\Vc^4_k)^{n^+_{lk}},\qquad\hfill
&f^{l+}_{\a^4\l}=f^{l+}_{\a^3\l}-n^+_{lk},\hfill\cr
x^{k+}_{\a^2\l}=(-1)^{n_{kl}^+}
x^{k+}_{\a^3\l}\otimes(q^{-1}\Vc_l^3/\Vc^2_l)^{n^+_{kl}},\hfill&
f^{k+}_{\a^2\l}=f^{k+}_{\a^3\l}-n^+_{kl},\hfill
\endmatrix$$
and the identity $n^+_{kl}+n^+_{lk}=-a_{lk}$.

\vskip3cm
{\eightpoint{
$$\matrix\format\l&\l&\l&\l\\
\phantom{.} & {\text{Michela Varagnolo}}\phantom{xxxxxxxxxxxxx} &
{\text{Eric Vasserot}}\\
\phantom{.}&{\text{D\'epartement de Math\'ematiques}}\phantom{xxxxxxxxxxxxx} &
{\text{D\'epartement de Math\'ematiques}}\\
\phantom{.}&{\text{Universit\'e de Cergy-Pontoise}}\phantom{xxxxxxxxxxxxx} &
{\text{Universit\'e de Cergy-Pontoise}}\\
\phantom{.}&{\text{2 Av. A. Chauvin}}\phantom{xxxxxxxxxxxxx} & 
{\text{2 Av. A. Chauvin}}\\
\phantom{.}&{\text{95302 Cergy-Pontoise Cedex}}\phantom{xxxxxxxxxxxxx} & 
{\text{95302 Cergy-Pontoise Cedex}}\\
\phantom{.}&{\text{France}}\phantom{xxxxxxxxxxxxx} & 
{\roman{France}}\\
&{\text{email: michela.varagnolo\@math.u-cergy.fr}}\phantom{xxxxxxxxxxxxx} &
{\text{email: eric.vasserot\@math.u-cergy.fr}}
\endmatrix$$
}}

\vskip1cm
\Refs
\widestnumber\key{ABC}

\ref\key{AK}\by Akasaka, T., Kashiwara, M.\paper 
Finite dimensional representations
of quantum affine algebras\jour Publ. RIMS\vol 33\yr 1997\pages 839-867\endref

\ref\key{B}\by Beck, J.\paper Braid group action and quantum affine algebras 
\jour Commun. Math. Phys.\vol 165\yr 1994\pages 555-568\endref

\ref\key{CG}\by Chriss, N., Ginzburg, V.\book Representation theory and complex
geometry\publ Birkh\"auser\publaddr Boston-Basel-Berlin\yr 1997\endref

\ref\key{CP}\by Chari, V., Pressley, A.\paper  Weyl modules for classical and quantum 
affine algebras\jour q-alg preprint\vol 0004174\endref

\ref\key{D}\by Damiani, I.\paper La R-matrice pour les alg\`ebres quantiques
de type affine non tordu\jour Ann. Sci. \'Ecole Norm. Sup. (4)
\vol 31\yr 1998\pages 493-523\endref

\ref\key{DD}\by Damiani, I., De Concini, C.\book Quantum groups and Poisson groups
(in Representations of Lie groups and quantum groups)
\ed Baldoni, Picardello\publ Longman Scientific and Technical\yr 1994\pages 1-45\endref

\ref\key{GRV}\by Ginzburg, V., Reshetikhin, N., Vasserot, E.\paper 
Quantum groups and flag varieties
\jour Contemp. math.\vol 175\yr 1994\pages 101-130\endref

\ref\key{GV}\by Ginzburg, V., Vasserot, E.\paper Langlands reciprocity for 
affine quantum groups of type $A_n$\jour Internat. Math. Res. Notices
\vol 3\yr 1993\pages 67-85\endref

\ref\key{K}\by Kashiwara, M.\paper Crystal bases of the modified quantized enveloping algebra 
\jour Duke Math. J.\vol 73\yr 1994\pages 383-413\endref

\ref\key{KT}\by Khoroshkin, S.M., Tolstoy, V.N.
\paper Twisting of quantum (super) algebras. Connection of Drinfeld's
and Cartan-Weyl realization for quantum affine algebras\jour hep-th
preprint\vol 9404036\endref

\ref\key{L1}\by Lusztig, G.\book Introduction to quantum groups 
\publ Birkh\"auser\publaddr Boston-Basel-Berlin \yr 1994\endref

\ref\key{L2}\by Lusztig, G.\paper Bases in equivariant K-theory
\jour Represent. Theory\vol 2\yr 1998\pages 298-369\endref

\ref\key{LS}\by Levendorskii, S.Z., Soibelman Ya.S.
\paper Some applications of quantum Weyl group I\jour J. Geom. Phys.
\vol 7\yr 1990\pages 241-254\endref

\ref\key{N1}\by Nakajima, H.\paper Instantons on ALE spaces, quiver varieties,
and Kac-Moody algebras\jour Duke Math. J.\vol 76\yr 1994\pages 365-415\endref

\ref\key{N2}\by Nakajima, H.\paper Quiver varieties and Kac-Moody algebras
\jour Duke Math. J.\vol 91\yr 1998\pages 515-560\endref

\ref\key{N3}\by Nakajima, H.\paper Quiver varieties and finite dimensional 
representations of quantum affine algebras \jour q-alg preprint
\vol 9912158\endref

\ref\key{V}\by Vasserot, E.\paper Affine quantum groups and equivariant 
$K$-theory\jour Transformation Groups\vol 3\yr 1998\pages 269-299\endref

\ref\key{VV}\by Varagnolo, M., Vasserot, E.\paper On the $K$-theory of
the cyclic quiver variety\jour Internat. Math. Res. Notices
\vol 18\yr 1999\pages 1005-1028\endref

\endRefs
\enddocument